\documentclass{amsart}

\newtheorem{theorem}{Theorem}[section]
\newtheorem*{maintheorem}{Theorem}
\newtheorem{lemma}[theorem]{Lemma}
\newtheorem{proposition}[theorem]{Proposition}
\newtheorem{corollary}[theorem]{Corollary}

\theoremstyle{definition}

\newtheorem{example}[theorem]{Example}

\newtheorem{remark}[theorem]{Remark}
\newtheorem*{acknowledgement}{Acknowledgement}

\theoremstyle{remark}

\DeclareFontFamily{U}{wncy}{}
\DeclareFontShape{U}{wncy}{m}{n}{<->wncyr10}{}
\DeclareSymbolFont{mcy}{U}{wncy}{m}{n}
\DeclareMathSymbol{\Sh}{\mathord}{mcy}{"58}

\newcommand{\tmu}{\widetilde{\mu}}

\newcommand\mylabel[1]{\label{#1}}

\newcommand{\ZZ}{\mathbb{Z}}

\newcommand{\FF}{\mathbb{F}}

\newcommand{\PP}{\mathbb{P}}
\renewcommand{\AA}{\mathbb{A}}
\newcommand{\GG}{\mathbb{G}}

\newcommand  {\shA}     {\mathcal{A}}

\newcommand  {\lieg}     {\mathfrak{g}}
\newcommand  {\lieh}     {\mathfrak{h}}

\newcommand  {\uAut}    {\underline{\operatorname{Aut}}}
\newcommand  {\uHom}    {\underline{\operatorname{Hom}}}
\newcommand  {\uEnd}    {\underline{\operatorname{End}}}
\newcommand  {\uIsom}   {\underline{\operatorname{Isom}}}

\newcommand  {\tlieg}   {\widetilde{\mathfrak{g}}}
\newcommand  {\tf}      {\widetilde{f}}

\newcommand  {\Aut}     {\operatorname{Aut}}

\renewcommand{\cong}    {\equiv}

\newcommand  {\Ext}     {\operatorname{Ext}}

\newcommand  {\fppf}    {{\text{fppf}}}
\newcommand  {\Fr}      {\operatorname{Fr}}

\newcommand  {\Gal}     {\operatorname{Gal}}
\newcommand  {\GL}      {\operatorname{GL}}

\newcommand  {\Hom}     {\operatorname{Hom}}

\newcommand  {\Ig}      {\operatorname{Ig}}

\newcommand  {\iso}     {\simeq}

\newcommand  {\Lie}     {\operatorname{Lie}}

\newcommand  {\lra}     {\longrightarrow}

\newcommand  {\maxid}   {\mathfrak{m}}

\renewcommand{\O}       {\mathcal{O}}
\newcommand  {\odd}     {{\operatorname{odd}}}

\newcommand  {\ord}     {{\operatorname{ord}}}

\newcommand{\p}         {\phantom{{}^*}}

\newcommand  {\Pic}     {\operatorname{Pic}}

\newcommand  {\quadand} {\quad\text{and}\quad}

\newcommand  {\ra}      {\rightarrow}

\newcommand  {\Sing}    {\operatorname{Sing}}
\newcommand  {\SL}      {\operatorname{SL}}

\newcommand  {\Spec}    {\operatorname{Spec}}

\def\mydate{\number\day\space\ifcase\month \or January\or February\or March\or 
April\or May\or June\or July\or
August\or September\or October\or November\or December\fi \space\number\year}

\DeclareFontFamily{U}{wncy}{}
\DeclareFontShape{U}{wncy}{m}{n}{<->wncyr10}{}
\DeclareSymbolFont{mcy}{U}{wncy}{m}{n}
\DeclareMathSymbol{\Sh}{\mathord}{mcy}{"58}

\usepackage{amscd,amssymb,amsmath}
\usepackage[matrix,arrow]{xy}

\begin{document}

\title[The N\'eron model over the Igusa curves]
      {The N\'eron model over the Igusa curves}

\author[Christian Liedtke]{Christian Liedtke}
\address{Mathematisches Institut, Heinrich-Heine-Universit\"at,
Universit\"atsstr. 1, 40225 D\"usseldorf, Germany}
\curraddr{Department of Mathematics, Stanford University,
450 Serra Mall, Stanford  CA 94305, USA} 
\email{liedtke@math.stanford.edu}

\author[Stefan Schr\"oer]{Stefan Schr\"oer}
\address{Mathematisches Institut, Heinrich-Heine-Universit\"at,
Universit\"atsstr. 1, 40225 D\"usseldorf, Germany}
\curraddr{}
\email{schroeer@math.uni-duesseldorf.de}


\dedicatory{Revised version, 3 May 2010}

\begin{abstract}
We analyze the geometry of rational $p$-division points in degenerating
families of elliptic curves in characteristic $p$.
We classify the possible Kodaira symbols and determine for the Igusa
moduli problem the reduction type of the universal curve.
Special attention is paid to characteristic $2$ and $3$, where
wild ramification and stacky phenomena show up.
\end{abstract}

\maketitle
\tableofcontents
\renewcommand{\labelenumi}{(\roman{enumi})}

\section*{Introduction}

Let $R$ be a discrete valuation ring of characteristic $p>0$, with function field $R\subset K$
and residue field $k=R/\maxid_R$,
which for simplicity we assume to be algebraically closed.
Suppose $E_K$ is an elliptic curve  containing
a \emph{rational $p$-division point} $z\in E_K$, that is a $K$-rational point of order $p$.
The goal of this paper is to analyze the N\'eron model $E\ra\Spec(R)$ and the geometry of the specialization
$\overline{\left\{z\right\}}\subset E$. 
The corresponding problem for rational division points of order prime to $p$ is classical and is
basically answered by the Ogg-Shafarevich Criterion.
For rational $p$-division points the situation turns out to be rather different.

Our point of departure are  the following   questions:
\begin{enumerate}
 \item When does there exist a rational $p$-division point $z\in E_K$ at all?
 \item What are the possible \emph{Kodaira symbols} describing the reduction types for $E$?
 \item Can the  specialization of $z$ into the closed fiber  $E_k$ be nonzero?
 \item  And for characteristic two and three, which are somewhat special for our considerations,
    can the class of $z$ in the \emph{component group} $\Phi_k$ be nonzero as well?
\end{enumerate}
Naturally these questions are closely related to the universal  family $U\ra\Ig(p)$ over the 
\emph{Igusa curve},
which, roughly speaking, parametrizes elliptic curves endowed with a rational $p$-division
point on the Frobenius pullback.

For the existence of a rational $p$-division point we analyze the kernel of the
multiplication-by-$p$ map.
For ordinary elliptic curves, this is a twisted form of $\mu_p\oplus(\ZZ/p\ZZ)$
and a rational $p$-division point exists if and only if this twisted form is trivial.
We classify twisted forms of this group scheme in terms of nonabelian cohomology
and relate this classification to the Hasse- and the $j$-invariant of an elliptic
curve.

In characteristic $p\geq5$
it turns out that the answer to the other questions depends on the congruence class 
of the characteristic modulo $12$.
Once we have an elliptic curve with
additive reduction and a rational $p$-division point with nonzero specialization, 
Frobenius pullbacks provide elliptic curves with the same reduction type as before and
zero specialization of arbitrary \emph{osculation number}.
Hence we are mainly interested in nonzero specialization of the $p$-division point and
the first main result is:

\begin{maintheorem}
Suppose $p\geq 5$.
An elliptic curve $E_K$ with additive reduction and containing a rational $p$-division point 
with nonzero specialization
exists precisely for the following reduction types:
\vspace{-0.8em}
$$
\begin{array}[t]{|l *{6}{|c}|}
\hline
\text{\rm congruence   mod $12$} & \text{\rm reduction type} \\
\hline
p\cong 1 & {\rm I}_0^*\\
\hline
p\cong 5 & {\rm II},  {\rm IV},  {\rm I}_0^*,{\rm IV}^*,  {\rm II}^*\\
\hline
p\cong 7 & {\rm III},  {\rm I}_0^*,  {\rm III}^*\\
\hline
p\cong 11 & {\rm II}, {\rm III},{\rm IV},  {\rm I}_0^*,{\rm IV}^*, {\rm III}^*,{\rm II}^*\\
\hline
\end{array}
$$
In this case $E_K$ has potentially supersingular reduction.
\end{maintheorem}

The same result holds if one demands that the rational point of order $p$ exists only
on the Frobenius pullback.
We actually construct explicit examples of such curves in terms of Weierstrass equations.
Note however, that we do not know the explicit coordinates for the rational $p$-division points. 
This existence of our examples has  strong consequences for  the universal elliptic curve over the Igusa curve:
Let $F$ be the function field of the Igusa curve $\Ig(p)$, and $U_F$ the corresponding universal
elliptic curve, and  $x\in\Ig(p)$ be a supersingular point.
Our second main result gives Weierstrass equations for $U_F$:

\begin{maintheorem}
Suppose $p\geq 5$. 
Then for  a suitable uniformizer $t\in\O_{\Ig(p),x}^\wedge$
the universal elliptic curve $U_F$ over the completion at $x$
is given by the following Weierstrass equations:
\vspace{-0.8em}
$$
\begin{array}[t] {|*{5}{c|}}
\hline
j(x) & p & \text{\rm Weierstrass equation} & U_F & U_F^{(p)}\\
\hline
0 & \cong -1 \!\!\!\mod 3 & y^2=x^3+t^{(p-5)/6}x + t^{-1} & {\rm III}^* & {\rm III}\\
\hline
1728 & \cong -1\!\!\!\mod 4 &y^2=x^3+t^{-1}x+t^{(p-7)/4} & {\rm II}^* & {\rm II}\\
\hline
\neq 0,1728 & \text{\rm all $p$} &y^2=x^3+at^{-2p}x+(b+t^{(p-1)/2})t^{-3p} & {\rm I}^*_0 & {\rm I}^*_0\\
\hline
\end{array}
$$
Here $a,b\in k$ are scalars so that the elliptic curve $y^2=x^3+ax+b$ has 
supersingular $j$-invariant $j=j(x)$.
\end{maintheorem}

Note that Ulmer \cite{Ulmer 1990a} gave Weierstrass equations with 
coefficients in $F$ for the universal curve $U_F$,
which rely on relations between Eisenstein series and are of somewhat
implicit nature. 
Our Weierstrass equations are explicit, but are defined only over
various completions.

The situation in characteristic $p=3$ and $p=2$ 
is more complicated and, in some sense, entirely different.
This comes from the fact that, besides the valuation of a minimal discriminant $\nu(\Delta)$,
there is an additional numerical invariant $\delta\geq 0$, the \emph{wild part of the conductor}.

\begin{maintheorem}
Let $p=3$.
For the Kodaira symbols ${\rm II}$, ${\rm II}^*$, ${\rm III}$, ${\rm III}^*$, ${\rm IV}$, ${\rm IV}^*$, ${\rm I}_0^*$, 
there is an elliptic curve $E_K$ containing a rational $3$-division point
with nonzero specialization in $E_k$ and the given reduction type.
For ${\rm IV}$ and ${\rm IV}^*$, there are such examples with nonzero specialization in $\Phi_k$, and
examples with zero  specialization in $\Phi_k$.
In any case, the curve has potentially supersingular reduction.
\end{maintheorem}

We show by example that the property of having a rational $3$-division point with nonzero class in $\Phi_k$
might even be preserved under base changes of arbitrarily large degree.
The universal elliptic curve over $\Ig(3)$ has already been determined by Ulmer \cite{Ulmer 1990a} and we reprove
this result in our setup.

In characteristic two, the Igusa moduli problem is not representable, such that
stacky phenomena show up.
Now, there is no restriction on the reduction type and the elliptic curve may 
have potentially ordinary and potentially multiplicative reduction.

\begin{maintheorem}
Let $p=2$.
For all additive Kodaira symbols, there is an elliptic curve $E_K$ containing a rational point of order two
with nonzero specialization in $E_k$ and having the given   reduction type.
For the Kodaira symbols ${\rm III}$, ${\rm III}^*$, ${\rm I}^*_l$, $l\geq 0$, there
are such examples where the specialization has nonzero class in $\Phi_k$, and examples
with zero class in $\Phi_k$. 
\end{maintheorem}

To prove the preceding results, we analyze the behavior of the wild part of the conductor under small field extensions,
and then apply Ogg's Formula $\nu(\Delta)=2+\delta+(m-1)$ to determine the reduction type.
It turns out that the numerical invariants have enough variation so that the Kodaira symbols listed above appear.

\medskip
The article is organized as follows:
In Section \ref{twisted mu} we classify twisted forms of $\mu_p$ and describe their $p$-Lie algebras.
In Section \ref{p torsion} we classify the twisted forms of $\mu_p\oplus(\ZZ/p\ZZ)$ in terms of nonabelian
cohomology.
This preparatory work is used in Section \ref{extension class} to describe the $p$-torsion subgroup scheme
of an ordinary elliptic curve.
In particular, we answer when an elliptic curve has a rational $p$-division point in this abstract setup
and relate this to the Hasse- and the $j$-invariant of the curve.
In Section \ref{models sections} we prove that an elliptic curve with additive reduction 
has potentially supersingular reduction, provided that it contains a rational $p$-division 
point with trivial specialization into the component group.
This result restricts the possible additive reduction types depending on the congruence class of
$p$ modulo $12$.
In Section \ref{osculation hasse} we introduce the notion of osculation number, that is, the
order of tangency of a rational $p$-division point with the zero section, and compute it
in terms of the Hasse invariant.
This allows us to decide when rational $p$-division points have nonzero specialization in the closed
fiber of the N\'eron model without computing the coordinates of these points explicitly.
In Section \ref{quadratic twists} we determine how the reduction type of an elliptic curve 
in characteristic $p\geq5$ changes under twisting.
This will be needed later on in the construction of examples.
In Section \ref{reduction frobenius} we determine how the reduction types 
in characteristic $p\geq5$ change under Frobenius pullbacks, which we need for the analysis
of the Igusa moduli problem.
In Section \ref{torsion twists} we construct elliptic curves with reduction of type
${\rm I}_0^*$ and having a rational $p$-division point with nonzero specialization in 
the special fiber.
These curves are obtained as quadratic twists of certain pullbacks of the versal deformation
of a given supersingular elliptic curve.
In Section \ref{decreasing osculation} we
start from versal deformations of supersingular elliptic curves with $j=0$ or $j=1728$
and construct elliptic curves having rational $p$-division points and nonzero specialization
in the special fiber for the remaining reduction types.
At this point we have shown that all possibilities determined in Section \ref{models sections}
do exist in characteristic $p\geq5$.
In Section \ref{igusa neron} we use our results to determine the degeneration behavior of
the universal elliptic curve over the Igusa moduli problem in characteristic $p\geq5$.
We even determine equations of the N\'eron model over the Igusa curves around
its supersingular points.
In Section \ref{delta one}, we specialize to characteristic $2$ and $3$, where we
analyze the Galois action on torsion points attached to elliptic curves whose wild 
part $\delta$ of the conductor is nontrivial yet as small as possible, namely $\delta=1$.
In Section \ref{igusa 3} we use these results to establish existence 
of elliptic curves in characteristic $3$ having a rational $3$-division point 
with nonzero specialization in the closed fiber for all possible reduction types.
Also, we determine the N\'eron model over the Igusa curve in characteristic $3$.
In Section \ref{igusa stack} we specialize to characteristic $2$ and introduce tautological 
families.
Since the Igusa moduli problem is not representable in characteristic $2$, these families
are in some sense the best replacement for the universal object.
We determine their reduction types and their behavior under Frobenius pullbacks.
In Section \ref{igusa 2} we construct the missing reduction types as pullbacks from
tautological families.
In Section \ref{semistable reduction 2} we classify reduction types in 
case we do not have potentially supersingular reduction.

\begin{acknowledgement}
  We thank Matthias~Sch\"utt and the referee for helpful comments.
\end{acknowledgement}

\section{Twisted forms of $\mu_p$ and their torsors}
\mylabel{twisted mu}

Let $S$ be a base scheme of characteristic $p>0$, endowed with the $\fppf$-topology.
Consider the finite diagonalizable group scheme $\mu_p=\Spec(\O_S[\ZZ/p\ZZ])$, whose values on $\Spec(A)\ra S$ are
$$
\label{mu values}
\mu_p(A)=\left\{x\in A\mid x^p=1\right\}.
$$
In this section we shall discuss twisted forms $\widetilde{\mu}_p$ of $\mu_p$
and the corresponding group $H^1(S,\widetilde{\mu}_p)$ of isomorphism classes of $\widetilde{\mu}_p$-torsors.
The former occur ``in nature'' as  the kernel of the relative Frobenius
on ordinary elliptic curves.
The results will be used in  the next section, which contains an analogous analysis for the group scheme $\mu_p\oplus\ZZ/p\ZZ$.
Throughout, we assume for simplicity that $S=\Spec(R)$ is affine, and that
$\Pic(R)=0$; for example, $R$ could be a field, a local ring, or a polynomial ring over a field.

Let $\shA=\underline{\Aut}(\mu_p)$ be the sheaf of automorphisms of $\mu_p$.
A twisted form $\tmu_p$ of $\mu_p$ determines an $\shA$-torsor
$\uIsom(\mu_p,\tmu_p)$. Conversely, an $\shA$-torsor $T$ yields
a twisted form $\tmu_p=T\wedge^\shA \mu_p=(T\times\mu_p)/\shA$. Here the quotient
is with respect to the diagonal action of $\shA$ on the product.
This establishes a canonical bijection between the   cohomology set
$H^1(S,\shA)$ and the set of isomorphism classes of twisted forms of $\mu_p$.
This correspondence has nothing in particular to do with $\mu_p$;
rather, it gives a general classification of twisted forms of sheaves 
(see \cite{Grothendieck 1955}, Chapter V,  or  \cite{Giraud 1971}, Chapter III, Section 2.3.
An exposition in the context of Galois cohomology can be found
in \cite{Serre 1994}, Chapter I, \S 5).

Let us write down the sheaf of automorphism $\shA$:
We have a canonical map 
$$
(\underline{\ZZ/p\ZZ})^\times\lra\shA,\quad \zeta\longmapsto(a\longmapsto a^\zeta),
$$
and it follows from 
\cite{SGA 3b}, Expos\'e VIII, Corollary 1.6 that this map is bijective.
Since $S$ is of characteristic $p>0$, we have
$\mu_{p-1}=(\underline{\ZZ/p\ZZ})^\times$ and so we may also
write this bijection as $\mu_{p-1}\ra\shA$. Using the Kummer sequence
$1\ra \mu_{p-1}\ra\GG_m\stackrel{p-1}{\ra}\GG_m\ra 1$ and
our assumption $\Pic(S)=0$, we deduce that the coboundary map
$$
R^\times/R^{\times (p-1)}\lra H^1(S,\mu_{p-1})=H^1(S,\shA)
$$
is bijective. 
In other words, the set of isomorphism classes of twisted forms of $\mu_p$
is the abelian group $R^\times/R^{\times(p-1)}$. 
We call  $\tau\in R^\times$ the \emph{twist parameter} for
the corresponding twisted form $\tmu_p$.
To write down these group schemes
explicitly, we first determine their $p$-Lie algebras:

Write $\mu_p=\Spec(R[T]/(T^p-1))$, and let $I\subset R[T]/(T^p-1)$ be the principle ideal
generated by $T-1$. Then $\Lie(\mu_p)$ equals $(I/I^2)^\vee$, which is a free $R$-module of rank one,
with basis $u\in\Lie(\mu_p)$ the residue class  of $T-1$. Using
$$
T^\zeta - 1 = (T-1)(T^{\zeta-1}+T^{\zeta-2}+\ldots +1)\cong
(T-1)\zeta \quad\text{modulo $I^2$},
$$
we see that the scalars $\zeta\in\mu_{p-1}$ act  on $\Lie(\mu_p)$ by scalar multiplication with $\zeta$.
A straightforward computation shows that the $p$-fold composition of the
derivation $(T-1)\frac{d}{d(T-1)}$ equals itself. In other words, the $p$-th power operation
is given by  $u^{[p]}=  u$.
We now view $H^1(S,\shA)=R^\times/R^{\times(p-1)}$ as the set
of isomorphism classes of twisted forms $\tlieg$ of the $p$-Lie algebra $\lieg=\Lie(\mu_p)$.

\begin{proposition}
\mylabel{twisted lie}
Let $\tau\in R^\times$. Then the corresponding twisted form $\tlieg$
is the 1-dimensional $p$-Lie algebra with   basis $\widetilde{u}\in\tlieg$
so that $\widetilde{u}^{[p]}=\tau^{-1} \widetilde{u}$.
\end{proposition}

\proof
By definition, the $p$-Lie algebra $\tlieg$ is the invariant submodule of the $R$-module $\lieg\otimes_R R[X]/(X^{p-1}-\tau)$,
where the action of $\zeta\in\mu_{p-1}$ is induced by  $X\mapsto\zeta X$ and $u\mapsto \zeta u$.
Clearly, $\widetilde{u}=u\otimes X^{-1}$ is invariant, and a basis of the twisted form $\tlieg$.
Its $p$-th power is $\widetilde{u}^{[p]}= u^{[p]}\otimes X^{-p}= \tau^{-1}\cdot u\otimes X^{-1} =\tau^{-1}\widetilde{u}$. 
\qed

\medskip
We finally regard $H^1(S,\shA)=R^\times/R^{\times(p-1)}$ as the set
of twisted forms $\tmu_p$ of the finite group scheme  $\mu_p$.

\begin{proposition}
\mylabel{twisted mu group scheme}
Let $\tau\in R^\times$. Then the corresponding twisted form $\tmu_p$
is the finite group scheme whose values on $R$-algebras $A$ 
are $\tmu_p(A)=\left\{a\in A\mid a^p=0\right\}$, with composition law
\begin{equation}
\label{oort-tate}
a\star b = a+b +\frac{1}{\tau}\sum_{i=1}^{p-1} \frac{a^ib^{p-i}}{i!(p-i)!}.
\end{equation}
\end{proposition}

\proof
First observe that the functor $G\mapsto\Lie(G)$ induces an equivalence
between the category of finite flat group schemes of height $\leq 1$
and the category of  $p$-Lie algebras whose underlying module is
projective of finite rank.
An  inverse functor is $\lieg\mapsto\Spec(U^{[p]}(\lieg)^\vee)$,
where  $U^{[p]}(\lieg)$ is the   the universal enveloping
algebra $U(\lieg)$ modulo the relations $x^p-x^{[p]}$, $x\in\lieg$.
Multiplication and comultiplication in the dual $U^{[p]}(\lieg)^\vee$
are induced by the diagonal $U^{[p]}(\lieg)\ra U^{[p]}(\lieg)\otimes U^{[p]}(\lieg)$
and the multiplication in $U^{[p]}(\lieg)$, respectively.

The $p$-Lie algebra of $\tmu_p$ is $\tlieg=Ru$ as in Proposition \ref{twisted lie},
and we merely have to spell out the general construction outlined
in the preceding paragraph for this special case.
Clearly, $1,u,\ldots,u^{p-1}\in U^{[p]}(\lieg)$ is an $R$-basis.
Let $f_0,\ldots,f_{p-1}\in U^{[p]}(\lieg)^\vee$ be the dual basis.
For $0\leq r,s,n\leq p-1$ we compute
$$
(f_rf_s)(u^n)
=\langle f_r\otimes f_s, (u\otimes 1+ 1\otimes u)^n\rangle
=\sum_{i=0}^n{n\choose i}f_r(u^i)f_s(u^{n-i})
={n\choose r}\delta_{s,n-r},
$$
where $\delta_{s,n-r}$ is a Kronecker Delta, and consequently
\begin{equation}
\label{multiplication table}
f_rf_s= 
\begin{cases}
{r+s\choose r} f_{r+s} &\text{if $r+s<p$},\\
0 &\text{else}.
\end{cases}
\end{equation}
Now set $f=f_1$.  Formula (\ref{multiplication table}) inductively gives $f^i=i!f_i$ for $0\leq i\leq p-1$ and $f^p=0$.
The upshot is that $U^{[p]}(\lieg)^\vee=R[f]/(f^p)$ as $R$-algebra.
It remains to compute the comultiplication map.
Now recall that $u^p=\tau^{-1}u$; this implies
$$
f_n(u^i\otimes u^j)=f_n(u^{i+j})=
\begin{cases}
\tau^{-1} &\text{if $i+j=n+p-1$},\\
1      &\text{if $i+j=n$},\\
0      &\text{else}.
\end{cases}
$$
for all $0\leq n, i,j\leq p-1$.
Putting things together, we infer that the comultiplication 
in $U^{[p]}(\lieg)^\vee$ is given by
\begin{equation}
\label{comultiplication}
f_n\longmapsto
\sum_{i+j=n}\frac{f^i\otimes f^j}{i!j!}
+
\frac{1}{\tau}\sum_{i+j=n+p-1}\frac{f^i\otimes f^j}{i!j!},
\end{equation}
where the summation indices satisfy $0\leq i,j\leq p-1$.
The special case $n=1$ now yields our assertions.
\qed

\begin{remark}
Formula (\ref{oort-tate}) is due to   Tate and Oort
\cite{Oort; Tate 1970}, page 9, who derived it, from a different perspective and 
in a   more general setting,
 with methods of representation theory.  They actually obtained
a classification of group schemes of order $p$ over
rather general base rings.
A discussion of their results is contained in
\cite{Shatz 1986}. For further generalizations, see \cite{Raynaud 1974}.
\end{remark}

\begin{remark}
Consider the special case $\tau=1$, such that $\tmu_p=\mu_p$.
At first glance, it seems strange that the
Formula (\ref{oort-tate}) gives a composition law $a\star b$  (on elements with $a^p=b^p=0$) 
that looks astonishingly different from the original composition
law $x\cdot y$ (on elements with $x^p=y^p=1$). Things clear up if one uses, instead of $f\in U^{[p]}(\lieg)^\vee$,
the truncated exponential
$$
e=f_0+f_1+\ldots+f_{p-1}=1+f+\frac{f^2}{2!}+\ldots+\frac{f^{p-1}}{(p-1)!}.
$$
Then $e^p=1$, and a direct computation using (\ref{comultiplication}) shows that
the  comultiplication indeed satisfies $e\mapsto e\otimes e$.
\end{remark}

Now fix a twisting parameter $\tau\in R^\times$, and let $\tmu_p$ be the corresponding twisted form of $\mu_p$.
We seek to understand the group $H^1(S,\tmu_p)$ of isomorphism classes of $\tmu_p$-torsors.
There seems to be no obvious relation to the group $H^1(S,\mu_p)$, because
the automorphisms of $\mu_p$ act via outer automorphisms,
compare \cite{Serre 1994}, Proposition 43.
In case $\tau\in R^{\times(p-1)}$, we have $\tmu_p\simeq\mu_p$, and the Kummer sequence
yields an isomorphism
\begin{equation}
\label{kummer p}
R^\times/R^{\times p}\lra H^1(S,\mu_p).
\end{equation}
In case $\tau\not\in R^{\times(p-1)}$, however, there is no embedding
of $\tmu_p$ into any iterated extension of the standard group schemes $\GG_a$  or $\GG_m$,
because there is no homomorphism of $p$-Lie algebras from $\tlieg=\Lie(\tmu_p)$
to $\Lie(\GG_a)$ or $\Lie(\GG_m)$. 
This destroys any hopes for an easy direct computations of $H^1(S,\tmu_p)$ such as (\ref{kummer p}).

However, there is an approach using Weil restriction.
Set $R'=R[X]/(X^{p-1}-\tau)$, such that $T=\Spec(R')$ is the $\mu_p$-torsor with
$\tmu_p=T\wedge^\shA\mu_p$.
Let $f:T\ra S$ be the structure morphism, and consider the \emph{Weil restriction}
$H=f_*(\mu_{p,T})$. This is a finite commutative group scheme on $S$ whose values on $R$-algebras $A$
are given by
\begin{equation}
\label{h functor}
H(A)=\left\{x\in(A\otimes_R R')^\times\mid x^p=1\right\}.
\end{equation}
To understand it, consider the twisted forms $H_i=T\wedge^\shA\mu_p$ of $\mu_p$,
where the action of $\zeta\in\mu_{p-1}=\shA$ on $\mu_p$ is given by $x\mapsto x^{\zeta^i}$, for $0\leq i\leq p-2$.
Clearly, $H_0=\mu_p$ and $H_1=\tmu_p$.

\begin{proposition}
\mylabel{sum decomposition}
There is a direct sum decomposition $H=H_0\oplus H_1\oplus\ldots\oplus H_{p-2}$.
\end{proposition}

\proof
It suffices to check this on the level of $p$-Lie algebras.
Set $\lieh=\Lie(H)$. Obviously,
$\lieh=\Lie(\mu_{p,T})=R'u$, viewed as an $R$-module.
Whence $u, Xu,\ldots, X^{p-2}u\in\lieh$ constitutes an $R$-basis, and $(X^iu)^{[p]}=X^{ip}u^{[p]}=\tau^i(X^iu)$.
Setting $\lieh_i=RX^iu$, we obtain a direct sum decomposition of $p$-Lie algebras
$\lieh=\lieh_0\oplus\lieh_1\oplus\ldots\oplus\lieh_{p-2}$.
Using Proposition \ref{twisted lie}, we infer that $\lieh_i=\Lie(H_i)$.
\qed

\medskip
This decomposition can also be viewed as an eigenspace decomposition:
The abelian Galois group $\mu_{p-1}$ acts functorially on the $\FF_p$-vector space $H(A)$ described in (\ref{h functor})  via its action on $R'$.
Whence the functor $H=H_{\chi^0}\oplus H_{\chi}\oplus\ldots \oplus H_{\chi^{p-2}}$ decomposes into eigenspaces, which are
indexed by the characters $\chi^i:\mu_{p-1}\ra\mu_{p-1}$, $0\leq i\leq p-1$. 
Here $\chi$ is the tautological character $\chi(\zeta)=\zeta$, and $\chi^i(\zeta)=\zeta^i$.

\begin{proposition}
\mylabel{eigenspace decomposition}
We have $H_i=H_{\chi^i}$ for all $0\leq i\leq p-2$.
\end{proposition}

\proof
We first look at the $p$-Lie subalgebras $\lieh_i\subset\lieh$.
Clearly, $\lieh_i=RX^iu$ is $\mu_{p-1}$-invariant, and $\zeta\in\mu_{p-1}$ acts
via multiplication by $\zeta^i=\chi^i(\zeta)$.
It follows that $H_i\subset H$ is $\mu_{p-1}$-invariant, and $\zeta\in\mu_{p-1}$ acts
via multiplication by $\zeta^i=\chi^i(\zeta)$. Whence $H_i\subset H_{\chi^i}$.
Since the inclusion
$$
H=H_0\oplus\ldots\oplus H_{p-2} \subset H_{\chi^0}\oplus\ldots\oplus H_{\chi^{p-2}}=H
$$
is an equality, we conclude $H_i=H_{\chi^i}$.
\qed

\medskip
The   decomposition is inherited to the cohomology groups of $H=f_*(\mu_{p,T})$.
This leads to the desired computation of $H^1(S,\tmu_p)$:

\begin{theorem}
\mylabel{eigenspace cardinality}
The  cohomology group $H^1(S,\tmu_p)$ is the $\chi$-eigenspace inside the cohomology group $H^1(S,f_*(\mu_{p,T}))=R'^\times/R'^{\times p}$ with respect to the Galois action of $\mu_{p-1}$, where $\chi:\mu_{p-1}\ra\mu_{p-1}$ is the tautological 
character $\chi(\zeta)=\zeta$.
\end{theorem}

\proof
We have $\tmu_p=H_1=H_\chi\subset H=f_*(\mu_{p,T})$, whence $H^1(S,\tmu_p)$ is the $\chi$-eigenspace
of $H^1(S,f_*(\mu_{p,T}))$. Since $f:T\ra S$ is finite and flat and $T\times_ST$ is a disjoint sum of copies of $T$,
it follows that $R^1f_*(G)=0$ for every abelian group scheme $G$ on $T$.
Now, the Leray--Serre spectral sequence shows that the canonical map $H^1(T,\mu_{p,T})\ra H^1(S,f_*(\mu_{p,T}))$ is bijective.
Finally, the Kummer sequence gives  $H^1(T,\mu_{p,T})=R'^\times/R'^{\times p}$.
\qed

\bigskip
Let us explicitly compute $H^1(S,\tmu_p)$ in the following special
case:
Suppose $D'$ is a normal noetherian domain of characteristic $p>0$ with function field $D'\subset F'$, endowed
with a faithful $\mu_{p-1}$-action.
We make the assumptions that $D'$ is factorial, and
that $(D')^{\times p}=(D')^\times$; this holds, for example,
for polynomial rings over perfect fields.
Let $F\subset F'$ be the field of $\mu_{p-1}$-invariants, 
and set $S=\Spec(F)$ and $T=\Spec(F')$, such that
the projection $f:T\ra  S$ is a $\mu_{p-1}$-torsor;
we denote by $\tmu_p$ the corresponding twisted form
of $\mu_p$.
Let  $I$ be the set of points of codimension one in $\Spec(D')$,
or equivalently the set of prime elements in $D'$ up to units.
Then $F'^\times/D'^\times$ is the free abelian group generated by $I$, and
the   action of $\mu_{p-1}$ on the ring
$D'$ induces a  permutation action on the set $I$.

\begin{proposition}
\mylabel{dimension cardinality}
Assumptions as above. Let $I_{\operatorname{free}\!}\subset I$
be the subset on which $\mu_{p-1}$ acts freely.
Then $H^1(S,\tmu_p)$ is an $\FF_p$-vector space
whose dimension equals the cardinality of the quotient set $I_{\operatorname{free}\!}/\mu_{p-1}$.
\end{proposition}

\proof
By Theorem \ref{eigenspace cardinality},  we may view  
$H^1(S,\tmu_p)$ as the $\chi$-eigenspace of the Galois module $H^1(T,\mu_{p,T})$,
where   $\chi:\mu_{p-1}\ra\mu_{p-1}$ is the tautological character.
The factoriality of $D'$ implies that the Galois module
$H^1(T,\mu_{p,T})=F'^{\times}/F'^{\times p}$ is the $\FF_p$-vector space with basis $I$.
In other words,  we have to 
compute the  $\chi$-eigenspace of $\FF_p[I]=\bigoplus_{i\in I}\FF_p$.
Using the algebra splitting of the group algebra $\FF_p[\mu_{p-1}]=\prod_{i=0}^{p-2}\FF_p$
into 1-dimensional  eigenspaces, we see that each free
orbit in $I$ contributes a 1-dimensional subspace to the $\chi$-eigenspace
of $\FF_p[I]$,
whereas each nonfree orbit in $I$ contributes only to eigenspaces for characters
$\chi^i\neq\chi$.
\qed

\begin{example}
Let $D'=k[X]$ be the polynomial ring over
a perfect field  $k$  of characteristic $p>0$,
on which $\zeta\in\mu_{p-1}$ acts via $X\mapsto\zeta X$.
Then $F'=k(X)$ and $F=k(X^{p-1})$, such that $\tau=X^{p-1}\in F$
is the twist parameter for the twisted form $\tmu_p$.
The set $I$   can be viewed
as the set of irreducible polynomials $h\in k[X]$ up to invertible scalars.
It is convenient to choose for every such $h$ with $h(0)\neq 0$ a representant
with $h(0)=1$. Such polynomials then factor as
$$
h(X)=(1-\alpha_1 X)(1-\alpha_2X)\ldots (1-\alpha_d X)
$$
with reciprocal roots $\alpha_1,\ldots,\alpha_d\in\Omega$ in some
algebraic closure $k\subset\Omega$.
The Galois action is then given by $h(\zeta X)=(1-\zeta\alpha_1 X)\ldots (1-\zeta\alpha_d X)$.
It is easy to see that such a polynomial yields a free Galois orbit in $I$
if and only if it is not contained in any of the subrings $k[X^i]\subset k[X]$, where $i>1$ ranges
over the divisors of $d-1$.
\end{example}

\section{Twisted forms of $p$-torsion in elliptic curves}
\mylabel{p torsion}

Let $S$ be a scheme of characteristic $p>0$,
endowed with the $\fppf$-topology.
Consider the  finite abelian group scheme $G=\mu_p\oplus\ZZ/p\ZZ$
over $S$. It is endowed with the \emph{Weil pairing}
$$
\Phi:G\times G\lra\mu_p,\quad
((\mu,i),(\nu,j))\longmapsto\mu^j/\nu^i,
$$
which is obviously bilinear, alternating, and nondegenerate.
We are interested in $(G,\Phi)$ because it or its twisted forms
naturally occur as the group scheme of $p$-torsion of ordinary elliptic curves.
The goal of this section is to determine
the set of isomorphism classes of \emph{twisted forms} $(\widetilde{G},\widetilde{\Phi})$ of $(G,\Phi)$.

This set can be viewed
as the set $H^1(S,\shA)$ of isomorphism classes of $\shA$-torsors, 
where $\shA=\uAut(G,\Phi)$. 
Our first task is to compute this sheaf of automorphism groups.
Let $U\ra S$ be a faithfully flat morphism of finite presentation.
Each local endomorphism of $G=\mu_p\oplus\ZZ/p\ZZ$ over $U$  can be written as a matrix
$$
\begin{pmatrix}
\zeta & \nu\\
  & \xi
\end{pmatrix}
\in\Gamma(U,\shA)
$$
with $\zeta,\xi,\nu$ local sections from $\uEnd(\mu_p)$, $\uEnd(\ZZ/p\ZZ)$,
$\uHom(\ZZ/p\ZZ,\mu_p)$, respectively.
There is no term below the diagonal because $\uHom(\mu_p,\ZZ/p\ZZ)=0$.
Using the canonical identifications 
$$
\underline{\ZZ/p\ZZ}=\uEnd(\mu_p),\quad
\underline{\ZZ/p\ZZ}=\uEnd(\ZZ/p\ZZ),\quadand
\uHom(\ZZ/p\ZZ,\mu_p)=\mu_p,
$$
we may view $\zeta,\xi:U\ra\ZZ/p\ZZ$ as locally constant functions,
and $\nu$ is an element from $\Gamma(U,\O_{U})$ with $\nu^p=1$.
The action on $G=\mu_p\oplus\ZZ/p\ZZ$ is given by
\begin{equation}
\label{group action}
\begin{pmatrix}
\zeta & \nu\\
  & \xi
\end{pmatrix}
\begin{pmatrix}
\mu\\
n
\end{pmatrix}
=
\begin{pmatrix}
\mu^\zeta\nu^n\\
\xi n
\end{pmatrix},
\end{equation}
and the composition law is
$$
\begin{pmatrix}
\zeta & \nu\\
  & \xi
\end{pmatrix}
\circ
\begin{pmatrix}
\zeta' & \nu'\\
  & \xi'
\end{pmatrix}
=
\begin{pmatrix}
\zeta\zeta' & \nu'^{\zeta}\nu^{\xi'}\\
&    \xi\xi'
\end{pmatrix}.
$$
Obviously, an endomorphism is an automorphism if and only if we have
$\zeta,\xi\in(\ZZ/p\ZZ)^\times=\mu_{p-1}$.
A straightforward argument shows that an automorphism
respects the Weil pairing $\Phi$ if and only if $\zeta\xi=1$.
Summing up, we have
\begin{equation}
\label{automorphism group}
\Gamma(U,\shA)=\left\{
\begin{pmatrix}
\zeta & \nu\\ &\zeta^{-1}
\end{pmatrix}
\mid
\text{$\zeta\in\Gamma(U,\mu_{p-1})$ and $\nu\in\Gamma(U,\mu_p)$}\right\}.
\end{equation}
We deduce that  $\shA$  sits inside an extension of groups
\begin{equation}
\label{extension}
1\lra\mu_p\lra\shA\lra\mu_{p-1}\lra 1,
\end{equation}
where the maps on the left and right are given by
$$
\nu\longmapsto \begin{pmatrix}1&\nu\\&1\end{pmatrix}
\quadand
\begin{pmatrix}\zeta&\nu\\&\zeta^{-1}\end{pmatrix}\longmapsto \zeta,
$$
respectively. The surjection $\shA\ra\mu_{p-1}$ has an obvious splitting given
by
$$
s:\mu_{p-1}\lra\shA,\quad 
\zeta\longmapsto \begin{pmatrix}\zeta&1\\&\zeta^{-1}\end{pmatrix} .
$$
Using this splitting, we may view
$\shA$ as a semidirect product $\shA=\mu_p\rtimes_\phi\mu_{p-1}$, 
for some homomorphism $\phi:\mu_{p-1}\ra\uAut(\mu_p)$.
 The  latter is given by $\zeta\mapsto(\mu\mapsto\mu^{\zeta^2})$, which follows from  the formula
for conjugation
\begin{equation}
\label{conjugation}
\begin{pmatrix}
\zeta & \nu\\
 &\zeta^{-1}
\end{pmatrix}
\cdot
\begin{pmatrix}
\xi&\mu \\
 &\xi^{-1}
\end{pmatrix}
\cdot
\begin{pmatrix}
\zeta & \nu\\
 &\zeta^{-1}
\end{pmatrix}^{-1}
=
\begin{pmatrix}
\xi &\mu^{\zeta^2}\nu^{\zeta(\xi^{-1}-\xi)} \\
 &\xi^{-1}
\end{pmatrix}.
\end{equation}
Recall that $\uAut(\mu_p)=\mu_{p-1}$, such that we may view $\phi$
as the map $\mu_p\ra\mu_p$, $\zeta\mapsto\zeta^2$.
Obviously, $\phi$ is trivial if and only if $p-1$, which is the order of the group $\mu_{p-1}$, divides $2$.
Hence:

\begin{proposition}
\mylabel{group commutative}
The   sheaf of groups $\shA$ is commutative if and only if $p=2$ or $p=3$.
In this case, we have $H^1(S,\shA)=H^1(S,\mu_p)\oplus H^1(S,\mu_{p-1})$.
\end{proposition}

For $p=2$ and $p=3$ it is thus easy to compute the cohomology group $H^1(S,\shA)$   with Kummer sequences.

From now on, we   assume that $p\geq 5$ and shall apply the theory of nonabelian cohomology to compute 
the cohomology set $H^1(S,\shA)$.
Care has to be taken because the extension in (\ref{extension})  is noncentral.
In any case, we have an exact sequence of pointed sets
$$
H^0(S,\shA)\lra H^0(S,\mu_{p-1})\lra H^1(S,\mu_p)\lra H^1(S,\shA) \lra H^1(S,\mu_{p-1}).
$$
The outer maps are surjective, because $\shA\ra\mu_{p-1}$ has   a section.
In other words:

\begin{proposition}
\mylabel{cohomology surjective}
The canonical map $H^1(S,\mu_p)\ra H^1(S,\shA)$ is injective,
the canonical map $H^1(S,\shA)\ra H^1(S,\mu_{p-1})$ is surjective,
and $H^1(S,\mu_p)$ is the fiber over the class of the
trivial torsor in $H^1(S,\mu_{p-1})$.
\end{proposition}

To understand the other fibers of the surjection $H^1(S,\shA)\ra H^1(S,\mu_{p-1})$,
it is necessary to twist the groups in (\ref{extension}).
Let $T$ be an $\shA$-torsor.
The sheaf of groups $\shA$ acts on itself by conjugation $a\mapsto(x\mapsto axa^{-1})$.
Whence we obtain a new sheaf of groups $\widetilde{\shA}=T\wedge^\shA \shA$, which is a twisted
form of $\shA$. The conjugation action leaves $\mu_p\subset\shA$ stable,
and is trivial on the quotient $\mu_{p-1}$.
Hence we  obtain a twisted form  $\widetilde{\mu}_p$ and an extension of groups
\begin{equation}
\label{twisted}
1\lra\widetilde{\mu}_p\lra\widetilde{\shA}\lra\mu_{p-1}\lra 1.
\end{equation}
It turns out that this extension does not necessarily split.
Note that in our situation the notions of schematically split 
and group-theoretically split   coincide:

\begin{proposition}
\mylabel{the splittings}
If the morphism of schemes $\widetilde{\shA}\ra\mu_{p-1}$ admits a section,
then there is also a section that is a homomorphism of group schemes.
In any case, there is at most one section that is a homomorphism.
\end{proposition}

\proof
Suppose there is a section  of schemes, and choose
a  generator $\zeta\in\mu_{p-1}$. Let $a\in\widetilde{\shA}(S)$ be the image of  $\zeta$ under the section.
Then $a^{p-1}\in\widetilde{\shA}(S)$ lies over $1\in\mu_{p-1}$, in other words, $a^{p-1}\in\tmu_p(S)$.
Since $p$ annihilates the group scheme $\tmu_p$, there is some $b\in\tmu_p(S)$ 
with $b^{1-p}=a^{p-1}$, namely $b=a^{p-1}$.
Replacing $a$ by $ba$ we obtain $a^{p-1}=1$. Whence $a$ defines a section that is also
a homomorphism of group schemes.

If $a,a'$ are two sections that are homomorphisms, then $a/a'$ defines a homomorphism
of group schemes $\mu_{p-1}\ra\tmu_p$. Such homomorphism must be trivial because 
 $p-1$ annihilates the domain of definition and $p$ annihilates the range. 
The uniqueness statement follows.
\qed

\medskip
Let $S'\ra S$ be the total space of the $\mu_{p-1}$-torsor induced from the $\shA$-torsor $T$ via
the homomorphism $\shA\ra\mu_{p-1}$.
According to Proposition \ref{cohomology surjective}, the pullback of $T$ along $S'\ra S$
is induced form a unique $\mu_{p,S'}$-torsor, which we call $T'$.

\begin{proposition}
\mylabel{twisted splits}
The   extension of groups in (\ref{twisted}) splits if and only if 
the $\mu_{p,S'}$-torsor $T'$ is trivial.
\end{proposition}

\proof
The condition is sufficient: Suppose that $T'$ has a section.
Our task is to see that the surjection  $\widetilde{\shA}\ra\mu_{p-1}$ has a section 
that is a homomorphism of groups.
We may check this after replacing $S$ by a finite Galois covering, because
if  such a section exist, it is unique by Proposition \ref{the splittings},
and whence descends.
Replacing $S$ by the total space $S'$ of the induced  $\mu_{p-1}$-torsor,
we may assume that the $\shA$-torsor $T$ is induced by some $\tmu_p$-torsor $T'$,
which is trivial by assumption. Now we are twisting with the trivial $\shA$-torsor $T$,
and the resulting group extension is obviously split.

The condition is also necessary: Suppose that $T'$ is nontrivial.
Let $\tf:\widetilde{\shA}\ra\mu_{p-1}$ be the canonical projection, and 
fix a generator $\xi\in\mu_{p-1}$. We shall show that the $\tmu_p$-torsor $\tf^{-1}(\xi)$ is nontrivial.
Making a base change as in the preceding paragraph, we may assume that the $\shA$-torsor $T$
is induced by the nontrivial $\mu_p$-torsor $T'$.
Let $f:\shA\ra\mu_p$ be the original projection. 
Then the fiber is $\tf^{-1}(\xi)=T'\wedge^{\mu_p}f^{-1}(\xi)$.
According to the formula for conjugation (\ref{conjugation}), the $\nu\in\mu_p$ act on $f^{-1}(\xi)$ via
$$
\begin{pmatrix}
\xi &\mu\\
    &\xi^{-1}
\end{pmatrix}
\longmapsto
\begin{pmatrix}
\xi &\mu\nu^n\\
    &\xi^{-1}
\end{pmatrix},
$$
where $n=\xi^{-1}-\xi $ is an element from $\ZZ/p\ZZ$. Since $p\geq 5$, we have $n\neq 0$, such that $\nu\mapsto\nu^n$ is
an automorphism of $\mu_p$. Set $m=1/n$, and let $T''$ be the   $\mu_p$-torsor obtained from $T'$ via pulling back
along the automorphism $\nu\mapsto\nu^m$.
Then $\tf^{-1}(\xi)$ is obtained from $f^{-1}(\xi)$  by twisting with the nontrivial
$\mu_p$-torsor $T''$ with respect to the multiplication action of $\mu_p$ on $f^{-1}(\xi)\simeq\mu_p$.
Consequently, $f^{-1}(\xi)\simeq T''$ does not admit a section over $S$.
\qed

\medskip
The preceding proof actually gives the following information:

\begin{proposition}
Suppose $S$ is connected, and that the $\mu_{p,S'}$-torsor $T'$ 
is nontrivial. Then the image of $H^0(S,\widetilde{\shA})\ra H^0(S,\mu_{p-1})$
is the subgroup $\left\{\pm 1\right\}$.
\end{proposition}

\proof
The condition that the fiber $\tf^{-1}(\xi)\subset\widetilde{\shA}$, $\xi\in\mu_{p-1}$ admits a section is equivalent to
the vanishing of $n=\xi^{-1}-\xi$, that is, $\xi=\pm 1$.
\qed

\medskip
We now have everything to compute the set of isomorphism classes of $\shA$-torsors:
Fix an $\shA$-torsor $ T$, and consider the exact
sequence (\ref{twisted}) obtained by twisting with $ T$ with respect
to the conjugation action. According to \cite{Grothendieck 1955}, Section 5.6, we have:

\begin{theorem}
\mylabel{cohomology set}
The pointed set of isomorphism classes of $\shA$-torsors  in $H^1(S,\shA)$
with the same image in $H^1(S,\mu_{p-1})$ as $ T$ is in canonical bijection
to the group $H^1(S,\widetilde{\mu}_p)$ modulo the permutation action of $H^0(S,\mu_{p-1})$ coming from
the twisted extension $1\ra\tmu_p\ra\widetilde{\shA}\ra\mu_{p-1}\ra 1$.
\end{theorem}

\begin{remark}
The permutation action of the subgroup $\left\{\pm 1\right\}\subset H^0(S,\mu_{p-1})$  on $H^1(S,\tmu_p)$ is trivial.
This follows from \cite{Grothendieck 1955}, Proposition 5.4.1, because
its conjugation action on $\tmu_p$ is trivial, and its image in $H^1(S,\tmu_p)$ under
the coboundary map is trivial.
\end{remark}

\section{The extension class}
\mylabel{extension class}

Let $K$ be a field of characteristic $p>0$ and $E_K$ an elliptic curve over $K$.
In order to decide when $E_K$ has a \emph{rational $p$-division point}, that is,
a $K$-rational point of order $p$, we shall analyze the multiplication-by-$p$ map.
We recall that the \emph{Frobenius pullback} $E_K^{(p)}$ is defined by the cartesian diagram
$$
\begin{CD}
E_K^{(p)} @>>> E_K\\
@VVV @VVV\\
\Spec(K) @>>F> \Spec(K).
\end{CD}
$$
We obtain a factorization
$$
\begin{xy}
\xymatrix{
E_K \ar[rr]^{p}\ar[dr]_{F} && E_K\\
& E_K^{(p)}\ar[ur]_{V}
}
\end{xy}
$$
of the multiplication-by-$p$ into the relative Frobenius followed by the 
\emph{Verschiebung}.
The kernels $\ker(F)$ and $\ker(V)$ are finite and flat group schemes
of order $p$ over $K$, where $\ker(F)$ is infinitesimal.
We recall that $E_K$ is \emph{ordinary} if 
$\ker(V)$ is \'etale, that is, a twisted form of $\ZZ/p\ZZ$.
In any case, the kernel $E_K[p]$ of the multiplication-by-$p$ map sits
inside a short exact sequence
\begin{equation}
\label{kernel of p}
1\lra\ker(F)\lra E_K[p]\lra\ker(V)\lra 1,
\end{equation}
and we conclude from this discussion:

\begin{proposition}
 There exists a rational $p$-division point on $E_K$ if and only
 if the following two conditions are satisfied:
 \begin{enumerate}
  \item the group scheme $\ker(V)$ is isomorphic to $\ZZ/p\ZZ$, and
  \item the extension (\ref{kernel of p}) splits.
 \end{enumerate}
\end{proposition}

As explained in \cite{Katz; Mazur 1985}, Section 2.8 there exists
a canonical pairing between the kernel of an isogeny and the kernel of its
dual isogeny.
This implies that $E_K[p]$ isomorphic to its own Cartier dual and hence
$\ker(F)$ is the Cartier dual of $\ker(V)$.
In particular, $\ker(V)$ is isomorphic to $\ZZ/p\ZZ$ if and only if
$\ker(F)$ is isomorphic to $\mu_p$.

We recall from \cite{Katz; Mazur 1985}, Section 12.4 that the
\emph{Hasse invariant} $h$ of $E_K$ is defined as the induced linear mapping on Lie algebras 
$$
h=\Lie(V):\Lie(E_K^{(p)})\ra\Lie(E_K).
$$
Using the identification $\Lie(E_K^{(p)})=\Lie(E_K)^{\otimes p}$, we may regard the Hasse invariant as an element
in the one-dimensional $K$-vector space $\Lie(E)^{\otimes(1-p)}$.
From this we derive more explicit invariants: Choose a basis $u\in\Lie(E_K)$, such
that $h=\lambda u^{\otimes(1-p)}$ for some scalar $\lambda\in K$, which is unique up to
$(p-1)$.st powers.

The Hasse invariant determines the $p$-Lie algebra $\lieg=\Lie(E_K)=Ku$ up to isomorphism via $u^{[p]}=\lambda u$,
and consequently $\ker(F)=\Spec(U^{[p]}(\lieg)^\vee)$.
Clearly, $E_K$ is ordinary if and only if $\lambda\neq0$.
In turn, the Cartier dual is given  as a scheme by
$$
\ker(V)= \uHom(\ker(F), \GG_m)=\Spec k[u]/(u^p-\lambda u).
$$
Applying Proposition \ref{twisted lie} we conclude

\begin{proposition}
 \label{twisted mu hasse}
 If $\lambda\neq0$, then $\ker(F)$ is the twisted form of $\mu_p$ corresponding to the
 twist parameter $\lambda^{-1}\in K$.
 In particular, the group scheme $\ker(V)$ is isomorphic to $\ZZ/p\ZZ$ if and only if
 $\lambda$ lies in $K^{\times(p-1)}$.
\end{proposition}

\begin{proposition}
 \mylabel{splitting with j}
 Let $E_K$ be an ordinary elliptic curve over $K$.
 Then the following are equivalent:
 \begin{enumerate}
  \item $j(E_K)\in K^p$,
  \item there exists an elliptic curve $X_K$ over $K$, so that
    $X_K^{(p)}\iso E_K$, and
  \item the extension (\ref{kernel of p}) splits.
 \end{enumerate}
 For supersingular elliptic curves (i) and (ii) are always true, whereas (iii)
 never holds.
\end{proposition}

\proof
Let $E_K$ be ordinary.
To prove $(i)\Rightarrow (iii)$, choose an elliptic curve $Y_K$ over $K$ with
$j(Y_K)^p=j(E_K)$.
The separable isogeny $V:Y_K^{(p)}\ra Y_K$ shows that the extension 
(\ref{kernel of p}) splits for $Y_K^{(p)}$ and the splitting is given by 
the subgroup scheme $G=\ker(V)$.
To proceed, note that $E_K$ is a twisted form of $Y_K^{(p)}$.
In order to establish existence of an \'etale subgroup scheme of $E_K[p]$ it
thus suffices to check that $G$ is invariant under the automorphism group scheme 
of $Y_K^{(p)}$.
To check the latter, we may assume that $K$ is algebraically closed.
Now the invariance is clear because the automorphism group scheme is reduced
and $G$ is the reduction of the $p$-torsion subgroup scheme.

To prove $(iii)\Rightarrow (ii)$, let $G_K\subset E_K[p]$ be the subgroup scheme
defining the splitting of (\ref{kernel of p}), and set $X_K=E_K/G_K$.
Consider the following commutative diagram
$$
\begin{xy}
\xymatrix{
E_K \ar[rr]^{\times p}\ar[d]_{{\rm pr}} && E_K\\
X_K \ar[rr]^{F} \ar@{-->}[urr] && X_K^{(p)} \ar@{-->}[u]
}
\end{xy}
$$
The diagonal dotted arrow exists because $G_K\subset E_K[p]$ and the
vertical dotted arrow exists because $X_K\ra E_K$ is purely
inseparable.
Then $X_K{(p)}\ra E_K$ has degree one, hence is an isomorphism.
The implication $(ii)\Rightarrow (i)$ is trivial.

Now let $E_K$ be supersingular.
Then (\ref{kernel of p}) never splits because otherwise the embedding
dimension of $E_K[p]$ would be too large.
As explained in the proof of \cite{Katz; Mazur 1985}, Theorem 12.4.3,
the multiplication by $p$-map induces a canonical
isomorphism $E_K\iso E_K^{(p^2)}$ from which
$(i)$ and $(ii)$ follow immediately.
\qed

\medskip
We have seen that $E_K[p]$ is a twisted form of $G=\mu_p\oplus(\ZZ/p\ZZ)$
over $S=\Spec K$ that respects the Weil pairing $\Phi$.
As explained in Section \ref{p torsion}, $E_K[p]$ defines
an $\shA=\uAut(G,\Phi)$-torsor.
By Proposition \ref{cohomology surjective} there is a surjective
homomorphism of pointed sets
$$
H^1(S,\shA)\lra H^1(S,\mu_{p-1})\lra 1
$$
mapping the class of $E_K[p]$ to the class of $\ker(F)$.
This latter cohomology group and the class of $\ker(F)$ 
have been analyzed in Section \ref{twisted mu}
and their relation to the Hasse invariant is described
in Proposition \ref{twisted mu hasse}.
To determine the fiber over the class of $\ker(F)$
we proceed as in Section \ref{p torsion} and 
consider the $\shA$-torsor 
$T=\ker F\oplus \uHom(\ker(F),\GG_m)$.
The twisted form 
$\widetilde{\shA}=T\wedge^\shA\shA$
of $\shA$ is an extension of $\mu_{p-1}$ by a twisted form
$\widetilde{\mu}_p$ of $\mu_p$ as in
(\ref{twisted}).
For this specific choice of $\widetilde{\shA}$
the general machinery developed in
Section \ref{p torsion} simplifies:

\begin{theorem}
 The pointed set of isomorphism classes of $\shA$-torsors
 with image $\ker(F)$ in $H^1(S,\mu_{p-1})$ is in bijection
 with $H^1(S,\widetilde{\mu}_p)$.
 Moreover we have canonical isomorphisms of groups
 $$
   H^1(S,\widetilde{\mu}_p)\iso H^1(S,\uHom(\ker V, \ker F))\iso
   \Ext^1(\ker V, \ker F),
 $$
 identifying this pointed set with the group of twisted splittings
 of (\ref{kernel of p}), as well as the group classifying all
 extensions of $\ker(V)$ by $\ker(F)$. 
\end{theorem}

\proof
Let $S'\ra S$ be the total space of the $\mu_{p-1}$-torsor induced
from the $\shA$-torsor $T$ via the homomorphism
$\shA\ra\mu_{p-1}$ as in Section \ref{p torsion}.
The pullback of $T$ along $S'\ra S$ yields 
$\mu_{p,S'}\oplus(\ZZ/p\ZZ)_{S'}$, which is
induced from the trivial $\mu_{p,S'}$-torsor.
By Proposition \ref{twisted splits}
the sequence (\ref{twisted}) for $\widetilde{\shA}$ is
split.
In particular, the map from
$H^1(S,\widetilde{\mu}_p)$ to $H^1(S,\widetilde{\shA})$ 
is injective.

By construction, we have $\widetilde{\mu}_p=\uHom(\ker V,\ker F)$,
so that $H^1(S,\widetilde{\mu}_p)$ classifies twists of splittings
of (\ref{kernel of p}).
The identification of the group of twisted splittings with
the group of all extensions follows from the discussion
in \cite{Demazure; Gabriel 1970}, Chapter III, \S6.3.5
and 
\cite{Demazure; Gabriel 1970}, Chapter III, \S6, Corollaire 4.9.
\qed

\medskip
We stress that this result is due to the specific choice
of the $\shA$-torsor $T$.
In this case, the distinguished element of $H^1(S,\widetilde{\shA})$
corresponds to the split extension of $\ker(F)$ by $\ker(V)$.
In particular, the class of $E_K[p]$ in $H^1(S,\widetilde{\shA})$ equals 
this distinguished element if and only if $j(E_K)\in K^p$
thanks to Proposition \ref{splitting with j}.
Moreover, in the proof we have seen that
if $f:S'\ra S$ trivializes the $\mu_{p-1}$-torsor $\ker(F)$ then
$\widetilde{\mu}_p$ becomes isomorphic to $\mu_p$.
Hence we obtain $\widetilde{\mu}_p$ as a subgroup scheme of
the Weil restriction $f_\ast(\mu_{p,S'})$ and 
Theorem \ref{eigenspace cardinality} applies.

\section{N\'eron models and sections of order $p$}
\mylabel{models sections}

Throughout, we shall work in the following set-up:
Let $R$ be a henselian discrete valuation ring of characteristic $p>0$,
whose residue field $k=R/\maxid_R$ is algebraically closed, with field of fraction $R\subset K$.
Let us also fix a uniformizer $t\in R$.
Given an elliptic curve $E_K$ over $K$, we denote by $E\ra\Spec(R)$ its \emph{N\'eron model}, and by $E_k\subset E$ the  closed fiber.
Let   $\Phi_k=E_k/E_k^0$ be the group of connected components of the closed fiber $E_k$.
Note that if $E_K$ has additive reduction, then the possible orders for $\Phi_k$ are  $1,2,3,4$.
We refer to \cite{Bosch; Luetkebohmert; Raynaud 1990} as general reference for the theory
of N\'eron models.

Suppose there is a rational $p$-division point $z\in E_K$.
Let $G_K\subset E_K$ be the subgroup scheme generated by $z$,
and consider its schematic closure $G\subset E$.
 
\begin{lemma}
\mylabel{subgroup scheme}
The subscheme $G\subset E$ is a subgroup scheme, and the
structure morphism $G\ra\Spec(R)$ is  flat and finite of degree  $p$.
\end{lemma}

\proof
Clearly, $G$ is reduced and the structure morphism $G\ra\Spec(R)$ is flat.
The N\'eron mapping property yields a morphism of group schemes $\varphi:(\ZZ/p\ZZ)\ra E$
with $1_K\mapsto z$. Hence $G$ is the schematic image of $\varphi$, whence finite because $E\ra\Spec(R)$ is separated.
Its degree must be $p$, because the generic fiber has length $p$.

To see that $G\subset E$ is a subgroup scheme, it suffices to check that the multiplication map
$\mu:G\times G\ra E$ factors through $G\subset E$, according to  \cite{Shatz 1970}.
Since $G\times G\ra\Spec(R)$ is flat and finite, the inclusion $G_K\times G_K\subset G\times G$ is dense.
Since $E\ra\Spec(R)$ is separated, the multiplication map $\mu:G\times G\ra E$ is finite, and in particular closed.
Whence
$$
\mu(G\times G)=\mu(\overline{G_K\times G_K}) = \overline{\mu(G_K\times G_K)} = \overline{G_K} = G,
$$
such that $\mu$ factors through $G\subset E$ set-theoretically. Using that $G\times G$ is reduced,
we conclude that the schematic image $\mu(G\times G)\subset E$ is reduced as well,
and infer that $\mu$ factors through $G\subset E$ scheme-theoretically.
\qed

\medskip
In the preceding situation, it is   convenient to consider the Cartier dual $H=\underline{\Hom}(G,\GG_m)$,
which is actually easier to describe than $G$. Note that the   group scheme $H\ra\Spec(R)$ is finite flat of degree $p$, 
and recall that $t\in R$ denotes a uniformizer.

\begin{proposition}
\mylabel{cartier dual}
Both fibers of $H\ra\Spec(R)$ are infinitesimal group schemes. If $G_k$ is connected,
then $G_k\simeq H_k\simeq\alpha_p$. 
Moreover, the Lie algebra $\lieh=\Lie(H)$ is a free $R$-module of rank one,
and admits a basis $b\in\lieh$ satisfying $b^{[p]}=t^nb$ for some 
integer $n\geq 0$.
\end{proposition}

\proof
By construction, the generic fiber is $H_K=\mu_{p,K}$.
Since $H_K\subset H$ is dense, it follows that the closed fiber $H_k$ is connected.
Over the algebraically closed field $k$, there are only two connected group schemes of length $p$,
namely $\alpha_p$ and $\mu_p$. Only the former has a connected Cartier dual.
So if $G_k$ is connected, we must have $H_k\simeq\alpha_{p,k}$.

The Lie algebra $\lieh$ is a free module of rank one, because the fibers of $H\ra\Spec(R)$ are
infinitesimal  of length $p$. Choose an arbitrary basis $b\in \lieh$.
Then $b^{[p]}=fb$ for some $f\in R$, which is nonzero because $H_K=\mu_{p,K}$.
Write $f=t^ng$ for some unit $g\in R$. 
Since $R$ is strictly henselian, there exists an $h\in R$ with $h^{p-1}=g$.
Replacing $b$ by $h^{-1}b$, we find the desired basis. 
\qed

\medskip
Now back to our elliptic curve  $E_K$ and its N\'eron model $E\ra\Spec(R)$.
If $K\subset K'$ is a finite field extension, we denote by $R'\subset K'$ the integral
closure of $R\subset K'$. Then $R'$ is a henselian discrete valuation ring with
field of fraction $R'\subset K'$ and algebraically closed residue field $k= R/\maxid_R=R'/\maxid_{R'}$.
We shall denote by $E_{K'}=E_K\otimes_KK'$ the  induced elliptic curve over $K'$,
and by $E'\ra\Spec(R')$ its N\'eron model. Note that the canonical map $E\otimes_R R'\ra E'$ coming
from the N\'eron mapping property is, in general, not an isomorphism.

The preceding Proposition yields a  first restriction on N\'eron models in presence
of rational $p$-division points:

\begin{theorem}
\mylabel{potentially supersingular}
Let $E_K$ be an elliptic curve over $K$.
Suppose the Frobenius pullback $E_K^{(p)}$ contains a rational $p$-division point
whose class in $\Phi_k$ is zero, and that
$E_K$ has additive reduction.
Then $E_K$ has potentially supersingular reduction.
\end{theorem}

\proof
For characteristic $p\neq 2$, this easily follows from the representability of the Igusa moduli problem.
The following argument works in general:
Replacing $E_K$ by $E_K^{(p)}$, we may assume that the rational $p$-division point already lies on $E_K$.
Choose  a finite field extension $K\subset K'$  over which $E_{K'}$ acquires semistable reduction.
The N\'eron mapping property yields a morphism $f:E\otimes_R R'\ra E'$, which is the identity over $K'$.
Since all homomorphisms from $\GG_a$ into $\GG_m$ or elliptic curves are zero,
$f$ maps  the connected component of the closed fiber of $E\otimes_R R'$ to the origin.

Fix a rational $p$-division point $z\in E_K$ and let $S_{z}\subset E$ be its closure.
Then the schematic image $f(S_z\otimes_R R')\subset E'$ is a section inducing a point of order $p$ over $K'$ and passing through the origin of the closed fiber.
Let $G'\subset E'$ be the closed subscheme generated by this section.
Its generic fiber is cyclic of order $p$, whereas
the closed fiber is connected.
According to Proposition \ref{cartier dual}, the closed fiber is isomorphic to $\alpha_p$.
Now suppose that $E_K$ has either potentially ordinary or potentially multiplicative reduction.
Then the connected component of the origin in $E'_k$ would be an ordinary elliptic curve or $\GG_m$.
But these  group schemes do not contain $\alpha_p$, a contradiction.
\qed

\begin{remark}
  If a rational $p$-division point has non-zero specialization into $\Phi_k$,
  then $p$ divides the order of $\Phi_k$.
  In case of additive reduction $\Phi_k$ is of order at most $4$.
  In particular, the assumption of the theorem on $\Phi_k$ 
  is automatically fulfilled for $p\geq5$.
\end{remark}

\medskip
Using information from tables of reduction types 
(for example  in \cite{Silverman 1994}, Chapter IV, \S 9),
we obtain the following more specific consequences:

\begin{corollary}
\mylabel{restriction i}
Let $E_K$ be an elliptic curve over $K$. Suppose that $E_{K}^{(p)}$ contains a rational $p$-division point, and that $p\geq 3$.  Then  the reduction type of $E_K$
is not ${\rm I}^*_l$ with $l\geq 1$.
\end{corollary}

\proof
If the reduction type is ${\rm I}^*_l$ with $l\geq 1$, then the $j$-invariant of $E_K$ is not contained in $R$. Consequently $E_K$ has potentially multiplicative reduction, in contradiction to Theorem 
\ref{potentially supersingular}.
\qed

\begin{corollary}
\mylabel{restriction p}
Let $E_K$ be an elliptic curve over $K$. Suppose that $E_{K}^{(p)}$ contains a rational 
$p$-division point.
If $p\geq5$ and if the reduction type is ${\rm II}$, ${\rm IV}$, ${\rm IV}^*$ or ${\rm II}^*$,
then we have $p\cong -1$ modulo $3$.
If $p\geq3$ and if the reduction type is ${\rm III}$ or ${\rm III}^*$, then $p\cong -1$ modulo $4$.
\end{corollary}

\proof
Let $j_k\in R/\maxid_R$ be the  residue class of the $j$-invariant of $E_K$, which must be a supersingular $j$-value by Theorem \ref{potentially supersingular}.
If the reduction type is ${\rm III}$ or ${\rm III}^*$, then $j_k=1728$ by the tables of reduction type.
According to \cite{Silverman 1986}, Chapter V, Example 4.5 this $j$-value is supersingular if and only if $p\cong -1$ modulo $4$.
If the reduction type is ${\rm II}$, ${\rm IV}$, ${\rm IV}^*$ or ${\rm II}^*$, then $j_k=0$,
and this is supersingular if and only if $p\cong -1$ modulo $3$, according to loc.\ cit.\ Example 4.4.
\qed

\begin{corollary}
\mylabel{restriction pi}
Let $E_K$ be an elliptic curve over $K$. Suppose that $E_{K}^{(p)}$ contains a rational 
$p$-division point, that $E_K$ has additive reduction, and that $p\cong 1$ modulo $12$.
Then $E_K$ has reduction type ${\rm I}_0^*$.
\end{corollary}

\proof
In light of Corollary \ref{restriction i} and Corollary \ref{restriction p},
the only remaining possibility for an additive reduction type is ${\rm I}_0^*$.
\qed

\section{Osculation numbers and Hasse invariant}
\mylabel{osculation hasse}

Let $E_K$ be an elliptic curve and $E\ra\Spec(R)$ be its N\'eron model,
say given by a minimal Weierstrass equation 
$y^2+a_1xy+a_3y=x^3+a_2x^2+a_4x+a_6$.
The group $E(K)$ comes along with a \emph{decreasing filtration} defined as follows \cite{Tate 1975}:
The subgroup $E_1(K)\subset E(K)$ comprises those $z\in E$ with vanishing specialization in $E_k$. 
The  coordinates of such points $z=(\lambda,\mu)$ have valuations 
$\nu(\lambda)=-2m$ and $\nu(\mu)=-3m$, and one defines
$$
E_m(K)=\left\{z\in E_1(K)\mid \text{$z=0$ or $\nu(\lambda)\leq -2m$}\right\}.
$$
Let us check that  this is independent of the chosen Weierstrass equation: 
Given a rational point $z\in E_K$, we write $S_z=\overline{\left\{z\right\}}$ for its closure
in $E$.
For each nonzero $z\in E_1(K)$, the scheme $S_z\cap S_0$ is a local Artin scheme,
and it is convenient to call its length the \emph{osculation number} of $z\in E_0(K)$.

\begin{lemma}
\mylabel{filtration well-defined}
The point $z\in E_1(K)$ has osculation number $m$ if and only if  $z\in E_m(K) - E_{m+1}(K)$.
\end{lemma}

\proof
Suppose $z=(\lambda,\mu)$ has osculation number $m$ and $\nu(\lambda)=-2n$.
By \cite{Tate 1975}, Theorem 4.2, the fraction $-x/y$, viewed as a variable, yields a uniformizer along the zero section for $E$.
The intersection scheme  $S_z\cap S_0$ has length $m$, and  is defined by $R[z]/(z,z-\lambda/\mu)$, which has length $\nu(\lambda/\mu)=n$.
We conclude $n=m$, and the result follows.
\qed

\medskip
Now we are interested in the osculation number for rational $p$-division points.
It turns out that this is closely related to the Hasse invariant.
Suppose for simplicity that $E_K$ has good reduction and consider the
factorization of the multiplication-by-$p$ morphism of N\'eron models
$$
\begin{xy}
\xymatrix{
E \ar[rr]^{p}\ar[dr]_{F} && E\\
& E^{(p)}\ar[ur]_{V}
}
\end{xy}
$$
As in Section \ref{extension class},
the \emph{Hasse invariant} $h$ of $E$ is defined as the induced linear mapping on Lie algebras 
$$
h=\Lie(V):\Lie(E^{(p)})\ra\Lie(E).
$$
Using the identification $\Lie(E^{(p)})=\Lie(E)^{\otimes p}$, we may regard the Hasse invariant as an element
in the invertible $R$-module $ \Lie(E)^{\otimes(1-p)}$.
From this we derive more explicit invariants: Choose a basis $u\in\Lie(E)$, such
that $h=\omega u^{\otimes(1-p)}$ for some scalar $\omega\in R$.
Then the \emph{vanishing order} $\nu(h)=\nu(\omega)\geq 0$ and the \emph{residue class} 
$[h]=[\omega]\in R/R^{\times(p-1)}$ 
do not depend on the choice of the basis.

The Hasse invariant determines the $p$-Lie algebra $\lieg=\Lie(E)=Ru$ up to isomorphism via $u^{[p]}=\omega u$,
and consequently $\ker(F)=\Spec(U^{[p]}(\lieg)^\vee)$.
It turn, the Cartier dual is given  as a scheme by
$$
\ker(V)=\uHom(\ker(F),\GG_m)=\Spec k[u]/(u^p-\omega u).
$$
This leads to the following observation:

\begin{proposition}
\mylabel{osculation p-torsion}
Suppose that $E_K$ has good reduction and that 
the Hasse invariant of $E$ has vanishing order $\nu(h)=p-1$.
Then the Frobenius pullback contains a rational $p$-division point  $z\in E_K^{(p)}$
with osculation number one. 
Its coordinates $z=(\lambda,\mu)$ in the Weierstrass model
have valuations $\nu(\lambda)=-2$ and $\nu(\mu)=-3$.
\end{proposition}

\proof
Since $R$ is strictly henselian, we may represent the Hasse invariant $h$ by the scalar $\omega=t^{p-1}$ for some
uniformizer $t\in R$. Clearly, $z\in \ker(V)$, and we just saw $\ker(V)=\Spec k[u]/(u^p-t^{p-1}u)$. Now the decomposition 
$$
u^p-t^{p-1}u=u\prod_{\zeta\in\mu_{p-1}(R)} (u-\zeta t)
$$
shows that the intersection $S_z\cap S_0$ has length one, and the statement follows from Lemma \ref{filtration well-defined}.
\qed

\begin{remark}
There is actually an explicit formula, of a somewhat implicit nature,
for the $x$-coordinate of   $p$-division points discovered by
Gunji \cite{Gunji 1976}. 
\end{remark}

\section{Reduction types under quadratic twists}
\mylabel{quadratic twists}

In this section we shall analyze the behavior of reduction types
under quadratic twists. 
Let $E_K$ be an elliptic curve over $K$, and  $W\ra\Spec(R)$
be its Weierstrass model, that is, the relative cubic defined by a minimal Weierstrass equation.
Choose a separable quadratic field extension $K\subset K'$,
and let $R\subset R'$ be the corresponding extension of discrete valuation rings.
Then the group $\left\{\pm 1\right\}$ acts on   $W$ via the sign involution,
and on $R'$ via the Galois involution.
Now consider the diagonal action on the product
$$
W'=W\times_{\Spec(R)}\Spec(R')=W\otimes_R R',
$$
and let $Y=W'/\left\{\pm 1\right\}$ be the quotient, which exists as a scheme
because $W'$ carries an ample invertible sheaf.
All our actions are via $R$-morphisms, so $Y$ is an $R$-scheme.
The scheme $Y$ is normal, the morphism $Y\ra\Spec(R)$ is proper,
and the canonical map $R\ra H^0(Y,\O_Y)$ is bijective.

Taking quotients commutes with passing to invariant open subsets, so
$Y_K= E_{K'}/\left\{\pm 1\right\}$ is nothing but the \emph{quadratic twist} of $E_K$ with respect
to the field extension $K\subset K'$. More abstractly,
$Y_K=\Spec(K')\wedge^{\left\{\pm 1\right\}} E_K$, where
we view $\Spec(K')$ as a $\left\{\pm 1\right\}$-torsor.
To emphasize this aspect, we   write  $\widetilde{E}_K=Y_K$ for this elliptic curve.
Note that  $K\subset K'$ is unique up to isomorphism if $p\neq 2$, because
$R$ is assumed to be strictly henselian.

\begin{proposition}
\mylabel{quadratic good}
If $E_K$ has good reduction and $p\neq 2$, then the quadratic twist $\widetilde{E}_K$ has
reduction type ${\rm I}_0^*$.
\end{proposition}

\proof
First note that the Weierstrass model coincides with the N\'eron model $E$,
which is a relative elliptic curve.
The fixed schemes for the action on $E$ and $\Spec(R)$ are given by the
$2$-torsion scheme and the closed point, respectively. Whence
the   fixed points on $E'=E\otimes_RR'$ are the 2-torsion points in the closed fiber $E'_k$.
If $p\neq 2$, then there are four such fixed points, whose images
on $Y$ are four rational double points of type $A_1$, which comprise $\Sing(Y)$.
Let $X\ra Y$ be the minimal resolution of singularities. Since $Y$ is Gorenstein and
$Y_k$ is irreducible,
the relative canonical class $K_{X/R}$ is trivial.
It follows that $X$ is the regular model of $\widetilde{E}_K$. 
The closed fiber  $Y_k$  has multiplicity two, because
it is birational to the quotient of the double curve $E\otimes_R(R'/tR')$.
We infer that $\widetilde{E}_K$ has reduction type ${\rm I}_0^*$.
\qed

\medskip
It is not difficult to determine the behavior under twists for arbitrary 
reduction types without geometry, by merely using  Weierstrass equations and Ogg's Formula, 
at least if $p\neq 2,3$, which we assume for the rest of this section. 
Then, Ogg's formula 
\cite{Silverman 1994}, Chapter IV, Formula 11.1
tells us that 
$$
\nu(\Delta) = \varepsilon + (m-1),
$$
where $\nu(\Delta)$ is the valuation of a minimal discriminant, $m$ denotes
the number of irreducible components of the closed fiber of $E$
and $\varepsilon$ is equal to $0,1,2$, depending on whether $E_K$ has
good, multiplicative or additive reduction.

\begin{proposition}
\mylabel{quadratic twist}
For $p\geq5$ the reduction types of $E_K$ and its quadratic twists $\widetilde{E}_K$
are related as in the following table:
\vspace{-0.8em}
$$
\begin{array}[t]{|l|c|c|c|c|c|c|c|c|c|}
\hline
E_K             & \p{\rm I}_m &{\rm II}&{\rm III}&{\rm IV}&\p{\rm IV}^* &\p{\rm III}^*&\p{\rm II}^* &\p{\rm I}_{m}^*\\
\hline&&&&&&&&\\[-2ex]
\widetilde{E}_K & \p{\rm I}_{m}^*&\p{\rm IV}^*&\p{\rm III}^*&\p{\rm II}^*&{\rm II}&{\rm III}&{\rm IV}& \p{\rm I_{m}} \\
\hline
\end{array}
$$
\end{proposition} 

\proof
Choose a minimal Weierstrass equation
\begin{equation}
\label{start equation}
y^2=x^3+a_4x+a_6
\end{equation}
for $E_K$ with coefficients $a_4,a_6\in R$.
According to \cite{Silverman 1986}, Chapter X, \S 6, Proposition 5.4,
the quadratic twist $\widetilde{E}_K$ has Weierstrass equation
\begin{equation}
\label{first equation}
y^2=x^3+t^2a_4x+t^3a_6,
\end{equation}
whose discriminant is $t^6\Delta$. 
Let $j\in K$ be the $j$-invariant of $E_K$. We first consider the case $\nu(j)\geq 0$.
Suppose that $\nu(\Delta)=0,2,3,4$, such that the reduction type of
$E_K$ is ${\rm I}_0$, ${\rm II}$, ${\rm III}$, ${\rm IV}$, respectively.
Then the   Weierstrass equation (\ref{first equation}) is minimal, and Ogg's Formula
implies that the reduction type of the quadratic twist $\widetilde{E}_K$ is ${\rm I}_0^*$, ${\rm IV}^*$, ${\rm III}^*$, ${\rm II}^*$.
Now suppose that $ \nu(\Delta)=6,8,9,10$, such that $E_K$ has  reduction type ${\rm I}_0^*$, ${\rm IV}^*$, ${\rm III}^*$, ${\rm II}^*$, respectively.
According to Lemma \ref{valuation bound}, the change of coefficients $x=t^2x'$ and $y=t^3y'$ yields another
integral Weierstrass equation
\begin{equation}
\label{second equation}
y^2=x^3+t^{-2}a_4x+t^{-3}a_6,
\end{equation}
which has discriminant $\widetilde{\Delta}=t^{-6}\Delta$, and is therefore minimal.
Ogg's Formula implies that the quadratic twist has reduction type  ${\rm I}_0, {\rm II}$, ${\rm III}$, 
${\rm IV}$, respectively.

It remains to treat the case $\nu(j)<0$.
Suppose that $E_K$ has reduction type ${\rm I}_m$, $m\geq 1$. 
Then $\nu(j)=-m$, and $a_4,a_6\in R$ are invertible,
by the table after \cite{Bosch; Luetkebohmert; Raynaud 1990}, Chapter 1, Section 1.5, Lemma 4.
It follows from the Tate algorithm (\cite{Tate 1975}, see also \cite{Silverman 1994}, Chapter 4, Section 9)
that the Weierstrass equation 
(\ref{first equation}) is minimal and has reduction of type ${\rm I}_l^*$ for some
$l\geq0$.
Then, Ogg's formula yields $l=m$.
Conversely, if $E_K$ has reduction type ${\rm I}_m^*$, then the 
Weierstrass equation (\ref{second equation}) is minimal,
with invertible coefficients, and the quadratic twist has reduction type ${\rm I}_m$.
\qed

\medskip
If $j(E_K)=0$, then the automorphism scheme of $E_K$ is isomorphic to $\mu_6$,
and we may also perform a \emph{sextic twist} with respect to
the generator
$$
u\in H^1(K,\mu_6)=K^\times/K^{\times 6}.
$$
If $E_K$ has Weierstrass equation $y^2=x^3+a_6$, then the sextic twist
is given by the Weierstrass equation $y^2=x^3+ua_6$, by \cite{Silverman 1986}, Chapter X, \S 6, Proposition 5.4.
Using $u^2$ instead of the generator $u\in H^1(K,\mu_6)$, we obtain the \emph{cubic twist}, which is given by
the Weierstrass equation $y^2=x^3+u^2a_6$.

\begin{proposition}
Suppose that $j(E_K)=0$ and $p\geq5$. 
Then the reduction types of $E_K$ and its cubic and sextic twists are related as in the following table:
\vspace{-0.8em}
$$
\begin{array}[t]{|l*{6}{|c}|}
\hline
E_K                   & \p{\rm I}_0 & {\rm II} & {\rm IV}    & \p{\rm I}_{0}^* & \p{\rm IV}^* & \p{\rm II}^* \\
\hline&&&&&&\\[-2ex]
\text{\rm cubic twist } & {\rm IV}  & \p{\rm I}_0^* & \p{\rm IV}^*  & \p{\rm II}^* &\p{\rm I}_0 &{\rm II}\\
\hline&&&&&&\\[-2ex]
\text{\rm sextic twist}  & {\rm II} & {\rm IV} & \p{\rm I}_0^* & \p{\rm IV}^* & \p{\rm II}^* & {\rm I}_0\\
\hline
\end{array}
$$
\end{proposition} 

If $j(E_K)=1728$, then the automorphism scheme of $E_K$ is isomorphic to $\mu_4$,
and we may also perform a \emph{quartic twist} with respect to
the generator
$$
u\in H^1(K,\mu_4)=K^\times/K^{\times 4}.
$$
If $E_K$ has Weierstrass equation $y^2=x^3+a_4x$, then the quartic twist
is given by the Weierstrass equation $y^2=x^3+ua_4x$, by \cite{Silverman 1986}, Chapter X, \S 6, Proposition 5.4.

\begin{proposition}
Suppose that $j(E_K)=1728$ and $p\geq5$. 
Then the reduction types of $E_K$ and its quartic twist
are related as in the following table:
\vspace{-0.8em}
$$
\begin{array}[t]{|l*{4}{|c}|}
\hline
E_K                 & {\rm I}_0   & {\rm III}   & \p{\rm I}_{0}^* & \p{\rm III}^* \\
\hline&&&&\\[-2ex]
\text{\rm quartic twist} & {\rm III} & \p{\rm I}_0^* & \p{\rm III}^*   & \p{\rm I}_0\\
\hline
\end{array}
$$
\end{proposition} 

The proofs for the preceding two propositions are as for Proposition \ref{quadratic twist};
we leave the actual computations to the reader.
In the proof of Proposition \ref{quadratic twist}, we have used the following fact:

\begin{lemma}
\mylabel{valuation bound}
Suppose $y^2=x^3+a_4x+a_6$ is a minimal Weierstrass equation.
If $\nu(\Delta)\geq 6$, then $\nu(a_4)\geq 2$ and $\nu(a_6)\geq 3$.
\end{lemma}

\proof
Using the formulas
$\Delta = -16(4a_4^3+27a_6^2)$ and $j=1728 (4a_4)^3/\Delta$
for the discriminant and $j$-invariant, together with Ogg's Formula,
one obtains the following table, which gives $\nu(\Delta)$ and the reduction type
in dependence on the valuations of $a_4,a_6\in R$:
\vspace{-0.8em}
$$
\begin{array}[t]{|l *{6}{|c}|}
\hline
\nu(a_4)  & 0  &  0   & \geq 1 & 1 &\geq 1 & \geq 2  \\ 
\hline&&&&&&\\[-2ex] 
\nu(a_6)  & 0  & \geq 1 &  0 & \geq 2 & 1 & 2   \\
\hline&&&&&&\\[-2ex]  
\nu(\Delta) & -n  & 0 & 0 & 3  &2 & 4 \\
\hline&&&&&&\\[-2ex] 
\text{reduction type} &{\rm I}_n&{\rm I}_0&{\rm I}_0&{\rm III }&{\rm II}&{\rm IV}\\
\hline
\end{array}
$$
The statement follows from this table.
\qed

\section{Reduction type under Frobenius pullback}
\mylabel{reduction  frobenius}

Let $E_K$ be an elliptic curve over $K$, and consider its Frobenius
pullback $E_K^{(p)}$.
In this section we describe the reduction type of the Frobenius pullback 
in terms of the reduction type of the original elliptic curve.
Since an additive reduction type changes to a semistable reduction type only
after a nontrivial Galois extension by  \cite{Serre; Tate 1968}, Corollary 3 to Theorem 2,
the following fact holds:

\begin{lemma}
\mylabel{frobenius semistable}
The Frobenius pullback $E_K^{(p)}$ has semistable reduction if 
and only if $E_K$ has semistable reduction.
\end{lemma}

Let $j\in K$ be the $j$-invariant of $E$, such that  $j^p\in K$ is the $j$-invariant
of the Frobenius pullback.

\begin{proposition}
\mylabel{frobenius multiplicative}
If $E_K$ has reduction type ${\rm I}_m$  for some $m\geq 0$,
then the Frobenius pullback $E_K^{(p)}$ has reduction type
${\rm I}_{pm}$.
\end{proposition}

\proof
We have $\nu(j)=-m$ and $\nu(j^p)=-pm$.
Using Lemma \ref{frobenius semistable} we infer that the Frobenius pullback
has reduction type ${\rm I}_{pm}$.
\qed

\begin{proposition}
\mylabel{frobenius pullback}
If $E_K$ has reduction type ${\rm I}_m^*$  for some $m\geq 0$ and $p\geq 3$,
then the Frobenius pullback $E_K^{(p)}$ has reduction type
${\rm I}_{pm}^*$.
\end{proposition}

\proof
It is easy to see that quadratic twists commute with Frobenius pullbacks.
So for $p\geq5$ the statement follows from Proposition \ref{quadratic twist}
and Proposition \ref{frobenius multiplicative}.
Also, if $m\geq1$, we can argue for all $p\geq3$ as follows:
then $\nu(j)=-m$ and $\nu(j^p)=-pm$, and we infer that the Frobenius pullback
has reduction type ${\rm I}_{pm}^*$. 
It remains to treat the case $p=3$ and $m=0$.
Then $\nu(\Delta)=6$ and the Tate algorithm reveals that a minimal Weierstrass
equation exists with $\nu(a_1)\geq1$, $\nu(a_2)\geq1$, $\nu(a_3)\geq2$,
$\nu(a_4)\geq2$ and $\nu(a_6)\geq3$.
Obviously, this equation is no longer minimal after Frobenius pullback
and we infer that the discriminant of a minimal Weierstrass equation
has $\nu(\Delta^{(p)})=3\cdot 6-12=6$.
From \cite{Silverman 1994}, Chapter IV, \S9, Table 4.1. we get the
reduction type ${\rm I}_0^*$.
\qed

\medskip
As we shall see in Section \ref{igusa 2}, this does not hold true in characteristic two.

\begin{proposition}
\mylabel{frobenius table}
Suppose $p\geq 5$. Then the reduction type of $E_K$ is related to the reduction type of its
Frobenius pullback as described the following table, according to the congruence class
of $p$ modulo $12$:
\end{proposition}

\vspace{-1.5em}
$$
\begin{array}[t]{|l*{8}{|c}|}
\hline
E_K                   
& {\rm I}_m    & \p{\rm I}_m^*    & {\rm II}    & {\rm III} & {\rm IV} &  \p{\rm IV}^*    & \p{\rm III}^* & \p{\rm II}^* \\
\hline&&&&&&&&\\[-2ex]
\text{$E_K^{(p)}$ for $p\cong 1$} 
& \p{\rm I}_{pm} & \p{\rm I}_{pm}^* & {\rm II}    & {\rm III} & {\rm IV} &  \p{\rm IV}^*    & \p{\rm III}^* & \p{\rm II}^* \\
\hline&&&&&&&&\\[-2ex]
\text{$E_K^{(p)}$ for $p\cong 5$} 
& \p{\rm I}_{pm} & \p{\rm I}_{pm}^* & \p{\rm II}^*    & {\rm III} & \p{\rm IV}^* &  {\rm IV}    & \p{\rm III}^* & {\rm II}\\
\hline&&&&&&&&\\[-2ex]
\text{$E_K^{(p)}$ for $p\cong 7$} 
& \p{\rm I}_{pm} & \p{\rm I}_{pm}^* & {\rm II}    & \p{\rm III}^* & {\rm IV} &  \p{\rm IV}^*    & {\rm III} & \p{\rm II}^*\\
\hline&&&&&&&&\\[-2ex]
\text{$E_K^{(p)}$ for $p\cong 11$} 
& \p{\rm I}_{pm} & \p{\rm I}_{pm}^* & \p{\rm II}^*    & \p{\rm III}^* & \p{\rm IV}^* &  {\rm IV}   & {\rm III} & {\rm II}\\
\hline
\end{array}
$$

\proof
We already verified the first two columns of the table.
Suppose now that $E_K$ has reduction type  of the form ${\rm II}, {\rm III},\ldots,{\rm II}^*$.
Then the reduction type is entirely determined by $1\leq \nu(\Delta)\leq 11$
via Ogg's Formula,
and we have $\nu(\Delta^{(p)})\cong p\nu(\Delta)$ modulo $12 $ by Tate's Algorithm.
 Reducing modulo $3$ and $4$, we see that
the  possible congruence  classes for the prime $p$ are $1,5,7,11$. 
A direct computation now yields the entries of the table.
\qed

\begin{remark}
Let $T\in\left\{{\rm II},{\rm III},{\rm IV}\right\}$ be a Kodaira symbol.
The change on the   Kodaira symbols  
in passing from $E_K$ to $E_K^{(p)}$ is   not difficult to remember: If 
$p\cong 1$ modulo $12$,   nothing changes.
If $p\cong -1$, then $T\leftrightarrow T^*$.
If $p\cong 5$  then
$T\leftrightarrow T^*$ for the ``even'' Kodaira symbols, and nothing else changes.
If $p\cong -5$,   then $T\leftrightarrow T^*$ for the ``odd'' Kodaira symbols, and nothing else changes.
\end{remark}

\section{$p$-torsion under quadratic twists}
\mylabel{torsion twists}

We keep the notation as in the preceding section.
Let $E_K$ be an elliptic curve over $K$.
Choose a separable quadratic field extension $K\subset K'$ and let $\widetilde{E}_K$ be 
the corresponding quadratic twist.
Passing    to  the quadratic twist   may turn closed points into rational points.
This relies on a useful fact from Galois theory:

\begin{lemma}
\mylabel{diagonal action}
Let $F\subset L$ be a finite abelian field extension 
with Galois group $H$.
The quotient of $\Spec(L)\times_{\Spec(F)}\Spec(L)$ by the diagonal action of $H$ is
isomorphic to the disjoint sum of $|H|$ copies of $\Spec(F)$.
\end{lemma}

\proof
Set $A=\prod_{\sigma\in H}L$. Under the isomorphism 
$L\otimes_F L\ra A$, $x\otimes y\mapsto (x\sigma(y))_{\sigma}$,
the diagonal tensor product action of $H$ on $L\otimes_F L$ corresponds to the
diagonal product action on $A$ given by  $\tau\cdot(z_\sigma)_{\sigma}=(\tau(z_\sigma))_{\sigma}$.
The corresponding invariant ring is $A^H=\prod_{\sigma\in H}F$, which gives our statement.
\qed

\medskip
Now  suppose $z\in E_K$ is a closed point, so that the corresponding closed subscheme
$\Spec\kappa(z)\subset E_K$ is invariant under the sign involution, and
whose residue field extension $K\subset\kappa(z)$ is isomorphic to $K\subset K'$.
Set $E_{K'}=E_K\otimes K'$, and let
$$
r:E_{K'}\ra E_K\quadand q:E_{K'}\ra \widetilde{E}_K
$$ 
be the canonical morphisms.
Then, by the preceding Lemma, the closed subscheme $q(r^{-1}(z))\subset \widetilde{E}_K$ is the disjoint
union of two rational points $y_1,y_2\in\widetilde{E}_K$. This leads to the following result:

\begin{proposition}
\mylabel{quadratic rational}
Suppose there is a \'etale subgroup scheme $G_K\subset E_K$ of length $p$,
containing a closed point $z\in G_K$ so that the field extension $K\subset\kappa(z)$
is of degree two.
Then $\widetilde{E}_K$ contains a rational $p$-division point.
\end{proposition}

\proof
First note that $G_K$ is a twisted form of $\ZZ/p\ZZ$. 
For $p=2$, all such twisted forms are trivial, such that all points of $G_K$ are rational.
Our assumptions thus imply $p\geq 3$.
Then the two field extensions
$K\subset\kappa(z)$ and $K\subset K'$ are isomorphic, because $R$ is strictly henselian.
We now check that $z\in G_K$, viewed as a closed subscheme, is invariant under the sign involution.
Write $G_K$ as the spectrum of $K[T]/(T^p-\tau T)$ for some $\tau\in K^\times$.
Then the sign involution acts via $T\mapsto -T$.
The closed point $z\in G_K$ corresponds to an irreducible quadratic factor in $K[T]$ of
$
T^{p-1}-\tau=\prod_{\zeta\in\mu_{p-1}}(T-\zeta\alpha),
$
where $\alpha\in\Omega$ is a root of this polynomial. The quadratic factors are of the form 
$$
(T-\zeta\alpha)(T-\zeta'\alpha)=T^2-(\zeta+\zeta')\alpha T+\zeta\zeta'\alpha^2,
$$
whence $\zeta'=-\zeta$, and the quadratic factor is invariant under $T\mapsto -T$.
In light of the discussion preceding the Proposition, the image $q(r^{-1}(z))\subset\widetilde{E}_K$ is the
disjoint union of two rational points, which are necessarily of order $p$.
\qed

\medskip
It is easy to see that quadratic twisting is compatible with isogenies:
If $E_K\ra E_K^\flat$ is an isogeny, we obtain an isogeny
$\widetilde{E}_K\ra\widetilde{E}_K^\flat$ of the same  degree between the corresponding
quadratic twists.
We now apply this to the inseparable isogeny $F:E_K\ra E_K^{(p)}$ of degree $p$.
Clearly, the induced isogeny on quadratic twists is also inseparable, and it follows that
the quadratic twist of a Frobenius pullback is the Frobenius pullback of the quadratic twist.
We simply write $\widetilde{E}_K^{(p)}$ for this elliptic curve.

\begin{proposition}
\mylabel{quadratic hasse}
Suppose $p\geq 3$, that $E_K$ has good reduction, and that the Hasse invariant of $E\ra\Spec(R)$ has vanishing order $(p-1)/2$.
Then the quadratic twist $\widetilde{E}^{(p)}_K$ contains a rational $p$-division point.
Its specialization in the closed fiber of the N\'eron model is nonzero
and lies in the connected component of the origin, that is, its class in
$\Phi_k$ is zero.
\end{proposition}

\proof
Let $E^{(p)}$ be the N\'eron model of $E_K^{(p)}$.
According to Lemma \ref{subgroup scheme}, there is a subgroup scheme $G\subset E^{(p)}$ of order $p$
that is generically \'etale.
Let $t\in R$ be a uniformizer. Then $t^{(p-1)/2}\in R$ represents the Hasse invariant of $E$,
and $G$ is isomorphic to  the spectrum of $R[T]/(T^p-t^{(p-1)/2}T)$.
Now let $R'\subset K'$ be the integral closure of $R\subset K'$, and choose a uniformizer $t'\in R'$ with $t'^2=t$.
Using the decomposition
$$
T^p-t'^{p-1}T= T \prod_{\zeta\in\mu_{p-1}} (T-\zeta t'),
$$
we infer that  $G'=G\otimes_R R'$ decomposes into $p$ sections for $E'=E\otimes_R R'$,
which are invariant under the diagonal $\left\{\pm 1\right\}$-action
and intersect pairwise transversally in the fixed point $0\in E'_k$.
We conclude that $\widetilde{E}_K^{(p)}$ contains a rational $p$-division point.
It remains to determine its specialization behavior.

Consider the normal surface $Y=E^{(p)}_{R'}/\left\{\pm 1\right\}$, and let $q:E^{(p)}_{R'}\ra Y$ be the quotient map.
Let $X\ra Y$ be the blowing-up of the four rational double points of type $A_1$ on $Y$.
As explained in the proof for 
Proposition \ref{quadratic good}, the N\'eron model $\widetilde{E}$ of $\widetilde{E}_K$
is obtained from $X$ by removing the strict transform of the closed fiber $Y_k\subset Y$.
Choose a rational $p$-division point $z\in \widetilde{E}_K^{(p)}$ and consider the closures $S_{z},S_0\subset \widetilde{E}_K^{(p)}$.
If the Artin scheme $q(S_{z})\cap q(S_0)$ is of length one, the strict transforms of $q(S_{z})$ and $q(S_{0})$
on $X$ will be disjoint and must pass through the same irreducible component of the closed fiber.
So the following Lemma concludes the proof.
\qed

\begin{lemma}
\mylabel{local singularity}
Let $S=\Spec(A)$ be a regular local 2-dimensional scheme in characteristic $p\neq 2$
endowed with a $\left\{\pm 1\right\}$-action   whose only fixed point is the closed point.
Let $C_1,C_2\subset S$ be   invariant regular curves intersecting transversely,
and $q:S\ra S/\left\{\pm 1\right\}$ be the quotient map.
Then $q(C_1)\cap q(C_2)$ has length one.
\end{lemma}

\proof
Without loss of generality we may assume that  $A=k[[x,y]]$, and that the action is given by
$x\mapsto-x$ and $y\mapsto -y$ (see, for example, \cite{Schroeer 2004}, Lemma 5.4). 
The invariant ring is then $k[[x^2,xy,y^2]]$.
Set $D_i=q(C_i)$, and let $a\in S/\left\{\pm 1\right\}$ be the closed point.  Then
 $q^{-1}(a)$ has length three, and the projections $C_i\ra D_i$ have degree two.
If the integral curve $D_i$ were nonnormal, then $q^{-1}(a)\cap C_i$ would have length $\geq 4$. 
This is impossible, so the $D_i$ are regular.
Let $n\geq 1$ be the length of  $D_1\cap D_2$.
The Nakayama Lemma implies that $q^{-1}(D_1\cap D_2)$ has length $\leq 3n$.
On the other hand,  $q^{-1}(D_1\cap D_2)\cap(C_1\cup C_2)$ has length $4n-1$. 
This gives us the estimate $4n-1\leq 3n$, and consequently $n=1$.
\qed

\medskip
It remains to find a discrete valuation ring $R$ and an elliptic curve $E_K$ over the function field
$R\subset K$ whose N\'eron model $E\ra\Spec(R)$ meets
the assumptions of the   Proposition \ref{quadratic hasse}. 

\begin{theorem}
\mylabel{nonzero specialization}
Let $k$ be an algebraically closed field of characteristic $p\geq 3$, and
$j_k\in k$ be  a supersingular $j$-value. Then there is an elliptic curve
$E_K$ over the field $K=k(t)$ with reduction type ${\rm I}_0^*$
so that $E_K^{(p)}$ contains a rational $p$-division point whose specialization
in the closed fiber of the N\'eron model is nonzero and lies in the connected component
of the origin. Moreover, the  $j$-invariant of $E_K$ lies in $R$ and has residue class  $j_k$.
\end{theorem}

\proof
Let $V_k$ be the supersingular elliptic curve with the given $j$-invariant $j_k$.
Set $A=k[u]_{(u)}$ and choose 
a versal deformation $V\ra\Spec(A)$   of $V_k$.
According to Igusa's Theorem (\cite{Katz; Mazur 1985}, Theorem 12.4.3), the Hasse invariant of $V$ has vanishing order one.
Now set $R=k[t]_{(t)}$ with $t=u^{(p-1)/2}$.
The Hasse invariant of the induced family $V\otimes_R R'$ has vanishing order $(p-1)/2$.
Let $E$ be the quadratic twist of $V\otimes_R R'$. 
According to Proposition \ref{quadratic good} and Proposition \ref{quadratic hasse},
the elliptic curve $E$ has all desired properties.
\qed

\begin{remark}
 Here and the the sequel we are concerned with the \emph{existence} of rational $p$-division
 points. It might be interesting to compute their coordinates \emph{explicitly}.
\end{remark}

\section{Decreasing osculation numbers}
\mylabel{decreasing osculation}

In this section we develop a method to produce rational   $p$-division points
by passing from one Weierstrass equation to another that is defined over a smaller field.
The set-up is as follows: Let $R'$ be a henselian discrete valuation ring
in characteristic $p\geq 5$ with algebraically closed residue field $k=R'/\maxid_{R'}$
and field of fractions $R'\subset K'$.
We also fix a uniformizer $t'\in R'$.

Let $E_{K'}$ be an elliptic curve over $K'$ with good reduction and N\'eron model $E'\ra\Spec(R')$.
Choose a Weierstrass equation of the form
\begin{equation}
\label{old equation}
y'^2=x'^3 +  a'_4x' + a'_6
\end{equation}
with coefficients $a'_4,a'_6\in R'$, such that the discriminant $\Delta'\in R'$ is invertible.
Making the substitutions $x'=t'^{-2}x$ and $y'=t'^{-3}y$ over $K'$, we obtain a new  
Weierstrass equation 
\begin{equation}
\label{new equation}
y^2 =x^3  + a_4 x + a_6,
\end{equation}
for  $E_{K'}$. Note that its coefficients $a_i=t'^ia'_i$
remain integral, such that 
the new Weierstrass equation still defines a relative cubic $C'\ra\Spec(R')$.
This cubic, however, is not the  Weierstrass model of its generic fiber, because its discriminant  
$t'^{12}\Delta'$ is not invertible.
Now suppose there is a subring $R\subset R'$ with $ a_4,a_6\in R$
so that the extension $R\subset R'$ is  finite and separable.
Replacing $R$ by its normalization, we conclude that $R$ is
another henselian discrete valuation ring, and the residue field
extension $R/\maxid_R\subset R'/\maxid_{R'}$ is bijective.
Let $R\subset K$ be the field of fractions.

Our new Weierstrass equation (\ref{new equation})
defines a relative cubic $C\ra\Spec(R)$ with $C\otimes_R R'\simeq C'$. Let $E_K$ be its generic fiber,
such that $E_K\otimes_KK'=E_{K'}$, and $E\ra\Spec(R)$ be its N\'eron model.

\begin{proposition}
\mylabel{additive reduction}
Under the preceding assumptions, the degree $d=[K':K]$   satisfies 
the divisibility condition $d\mid 12$. If  furthermore $d\neq 1$, then 
$E_K$ has additive reduction 
and the relative cubic $C\ra\Spec(R)$ is its Weierstrass model.
If $d=6,4,3,2$, then the reduction type of $E_K$ is ${\rm II}$, ${\rm III}$, 
${\rm IV}$, ${\rm I}_0^*$, respectively.
\end{proposition}

\proof
Let $\Delta\in R$ be the discriminant for (\ref{new equation}).
Since $C\otimes _R R'=C'$ by construction, we have $d\nu(\Delta)=12$.
Now suppose that $d\neq 1$, such that $\nu(\Delta)<12$. 
By Tate's Algorithm, the Weierstrass equation (\ref{new equation})
must be minimal, such that $C\ra\Spec(R)$ is the Weierstrass model of
its generic fiber. Since $\Delta\in\maxid_R$, the elliptic curve
$E_K$ has bad reduction. Since $E_K$ has potentially good reduction,
the reduction type must be additive.
Ogg's Formula $\nu(\Delta)=2+(m-1)$ implies the statement
on the reduction types.
\qed

\medskip
We now examine the behavior of  rational $p$-division points in our construction:

\begin{proposition}
\mylabel{osculation decreases}
Under the preceding assumptions, suppose the field extension $K\subset K'$ has  degree $d>1$.
Assume furthermore that $E_{K'}$ contains a rational  $p$-division point
with osculation number one.
Then $E_K$ contains a rational   $p$-division point
whose specialization into $E_k$ is nonzero.
\end{proposition}

\proof
Choose a rational $p$-division point $z\in E_{K'}$, say with
coordinates $z=(\lambda,\mu)$ with $\lambda,\mu\in K'$.
According to Proposition \ref{filtration well-defined}, the coordinates have valuations
$\nu(\lambda)=-2$ and $\nu(\mu)=-3$.
Consequently $u'^2\lambda,u'^3\mu\in R'$,   and the closure 
$S_z\subset C'$ of $z\in E_{K'}$  in the relative cubic 
defined by the new Weierstrass equation (\ref{new equation}) is a section over $R'$ disjoint from
the zero section.
Since $u'^2\lambda$ is invertible, it is also disjoint from the  singularity in $C'$.

Suppose for the moment that the $j$-invariant $j\in K$ of the elliptic curve $E_K$
is a $p$-th power. According to Proposition \ref{splitting with j}, there is
an \'etale subgroup scheme $G_K\subset E_K$ of order $p$.
Let $A_K=G_K-0$ be the complement of the origin,
and $A\subset C$ be its closure in the Weierstrass model $C\ra\Spec(R)$ of $E_K$.
Since $d\geq 1$, this Weierstrass model is defined by the Weierstrass equation
(\ref{new equation}), according to Proposition \ref{additive reduction}.
We saw in the  preceding paragraph   that
$A\otimes_R R'\subset C'=C\otimes_R R'$ is disjoint from the zero section   and the singularity
in $C'$, so the same holds for $A\subset C$.
We infer that $A\cup\{0\}$ coincides with the closure $G\subset E$
of $G_K$ in the N\'eron model $E\ra\Spec(R)$, such that $G$ is a relative group scheme
whose closed fiber is reduced at the origin. But  $R$ is strictly henselian,
so $G=(\ZZ/p\ZZ)_R$. Restricting to the generic fiber yields
the desired rational $p$-division point on $E_K$.

It remains to verify that the $j$-invariant $j\in K$ of $E_K$
is a $p$-th power. By assumption, $E_{K'}$ contains a rational  $p$-division point,
whence $j\in K'$ is a $p$-th power by Proposition \ref{splitting with j}.
The following Lemma ensures that $j\in K$ is already a $p$-th power.
\qed

\begin{lemma}
\mylabel{separable extension}
Let $F\subset E$ be a   field extension in characteristic $p>0$.
If this extension is separable, then the inclusion $F^p\subset E^p\cap F$
is a bijection.
\end{lemma}

\proof
It suffices to show that the canonical map $F^\times/F^{\times p}\ra E^\times/E^{\times p}$ is injective.
Via the Kummer sequence, we may regard  this  map as $H^1(F,\mu_p)\ra H^1(E,\mu_p)$.
Let $T$ be a nontrivial $\mu_p$-torsor over $F$, such that $T$ is a reduced scheme.
Since $F\subset E$ is separable, the induced torsor  $T\otimes_FE$ remains a reduced scheme,
hence is a  nontrivial torsor. Consequently, the map in question is injective.
\qed

\begin{proposition}
\mylabel{congruent 2}
Suppose $p\cong -1$ modulo $3$. Fix an integer $1\leq n\leq 5$,
and let  $E_K$ be the elliptic curve over $K$ defined by the Weierstrass equation
\begin{equation}
\label{weierstrass mod 3}
y^2=x^3+t^{n(p-5)/6}x + t^{-n}.
\end{equation}
Then $E_K^{(p)}$ contains a rational $p$-division point
whose specialization in the closed fiber of the N\'eron model is nonzero.
The $j$-invariant of $E_K$ lies in $R$ and reduces to $0\in k$. The reduction types are given by the following table:
\vspace{-0.8em}
$$
\begin{array}[t]{|l|c|c|c|c|c|}
\hline
n    &  1 & 2 & 3 & 4 & 5\\
\hline&&&&&\\[-2ex]
E_K  &  \p{\rm II}^* & \p{\rm IV}^* & {\rm I}_0^* & {\rm IV} & {\rm II}\\
\hline&&&&&\\[-2ex]
E_K^{(p)}  &  {\rm II} & {\rm IV} & {\rm I}_0^* & \p{\rm IV}^* & \p{\rm II}^*\\
\hline
\end{array}
$$
\end{proposition}

\proof 
We first treat the case $n=1$.
Let $t'$ be an indeterminate, and consider the family of elliptic curves
$y^2=x^3+t'x+1$ over $R'=k[[t']]$, which is a versal deformation 
for the supersingular elliptic curve $y^2=x^3+1$ with $j$-invariant $j_k=0$.
The base change $t'\mapsto t'^{p-1}$ yields the family of elliptic curves
$y^2=x^3+t'^{p-1}x+1$, whose Hasse invariant has vanishing order $p-1$,
according to Igusa's Theorem (\cite{Katz; Mazur 1985}, Theorem 12.4.3).
Whence the Frobenius pullback $y^2=x^3+t'^{p(p-1)}x+1$ contains
a rational $p$-division point, which has osculation number one 
by Proposition \ref{osculation p-torsion}.
Making the substitution $x=t'^{-2}x'$, $y=t'^{-3}y'$, we obtain
the new Weierstrass equation $y^2=x^3+t'^{p(p-1)+4}x +t'^6$,
which defines a relative cubic over $k[[t']]$. Since $p\cong -1$ modulo $3$,
this cubic is already defined over the subring $k[[t]]\subset k[[t']]$,
where $t=t'^6$. According to Proposition \ref{osculation decreases},
the elliptic curve
$$
y^2=x^3+t^{(p(p-1)+4)/6}x + t
$$
over $K=k((t))$ contains a rational $p$-division point whose specialization
in the N\'eron model is nonzero.
Its reduction type is ${\rm II}$ by Proposition \ref{additive reduction}.

It remains to identify this elliptic curve as the Frobenius pullback
of $E_K$.
To do this, write $p=6d-1$ for some integer $d\geq 1$.
Making the substitution $x=t^{2d}x'$, $y=t^{3d}y'$, we obtain
the Weierstrass equation
$$
y^2=x^3+t^{(p(p-1)+4)/6-4d}x + t^{1-6d},
$$
which is indeed the Frobenius pullback of (\ref{weierstrass mod 3}) in the case $n=1$ at hand.
According to Proposition \ref{frobenius table}, the curve $E_K$ has reduction type
${\rm II}^*$.

The elliptic curves for $n>1$ are obtained from the elliptic curve with $n=1$
via the base change $t\mapsto t^n$, and their reduction types
easily follow from Ogg's Formula.
\qed

\begin{proposition}
\mylabel{congruent 3}
Suppose $p\cong -1$ modulo $4$. Fix an integer $1\leq n\leq 3$,
and let  $E_K$ be the elliptic curve over $K$ defined by the Weierstrass equation
\begin{equation}
\label{weierstrass mod 4}
y^2=x^3+t^{-n}x + t^{n(p-7)/4}.
\end{equation}
Then $E_K$ contains a rational $p$-division point whose specialization in the closed fiber of the N\'eron model is nonzero.
Its $j$-invariant lies in $R$ and reduces to $1728\in k$, and
the reduction type is given by the following table:
\vspace{-0.8em}
$$
\begin{array}[t]{|l|c|c|c|}
\hline
n   &  1 & 2 & 3\\
\hline&&&\\[-2ex]
E_K  &  \p{\rm III}^* &   {\rm I}_0^* & {\rm III} \\
\hline&&&\\[-2ex]
E_K^{(p)}  &  {\rm III} &   {\rm I}_0^* & \p{\rm III}^* \\
\hline
\end{array}
$$
\end{proposition}

\proof
This is analogous to the proof for Proposition \ref{congruent 2}.
Consider the family of elliptic curves
$y^2=x^3+x+t'$ over $R'=k[[t']]$, which is a versal deformation 
for the supersingular elliptic curve $y^2=x^3+x$ with $j$-invariant $j_k=1728$.
The base change $t'\mapsto t'^{p-1}$ yields the family of elliptic curve
$y^2=x^3+x+t'^{p-1}$, and its  Frobenius pullback $y^2=x^3+x+t'^{p(p-1)}$ contains
a rational $p$-division point with osculation number one.
Making the substitution $x=t'^{-2}x'$, $y=t'^{-3}y'$, we obtain
the new Weierstrass equation $y^2=x^3+t'^{4}x +t'^{p(p-1)+6}$,
which defines a relative cubic over $k[[t']]$. Since $p\cong -1$ modulo $4$,
this cubic is already defined over the subring $k[[t]]\subset k[[t']]$,
where $t=t'^4$. According to Proposition \ref{osculation decreases},
the elliptic curve
$$
y^2=x^3+tx + t^{(p(p-1)+6)/4}
$$
over $K=k((t))$ contains a rational $p$-division point whose specialization
in the N\'eron model is nonzero.
Its reduction type is ${\rm III}$   according to Proposition \ref{additive reduction}

Write $p=4d-1$ for some integer $d\geq 1$.
Making the substitution $x=t^{2d}x'$, $y=t^{3d}y'$, we obtain
the Weierstrass equation
$$
y^2=x^3+t^{1-4d}x + t^{(p(p-1)+6)/4-6d},
$$
which is   the Frobenius pullback of (\ref{weierstrass mod 4}) in the case $n=1$.
According to Proposition \ref{frobenius table}, the curve $E_K$ has reduction type
${\rm III}^*$.
The elliptic curves for $n>1$ are obtained from the elliptic curve with $n=1$
via the base change $t\mapsto t^n$, and Ogg's Formula gives the reduction type.
\qed

\medskip
Now  let $j_k\in \overline{\FF}_p\subset R$ be an arbitrary supersingular $j$-value.
Choose $a,b\in\FF_{p^2}$ so that $y^2=x^3+ax+b$ defines
a supersingular elliptic curve with   $j$-invariant $j_k$.

\begin{proposition}
\mylabel{congruent arbitrary}
Assumptions as above. Let  $E_K$ be the elliptic curve over $K$ defined by the Weierstrass equation
\begin{equation}
\label{generic supersingular}
y^2=x^3+at^{-2p}x + (b+t^{(p-1)/2})t^{-3p}.
\end{equation}
Then $E_K$ contains a rational $p$-division point whose specialization in the closed
fiber of the N\'eron model is nonzero.
Its reduction type is ${\rm I}^*_0$, and its  $j$-invariant lies in $R$ and reduces to $j_k\in k$. 
\end{proposition}

\proof
This is analogous to the proof for Proposition \ref{congruent 2}.
Consider the family of elliptic curves
$y^2=x^3+ax+b+t'$ over $R'=k[[t']]$, which is a versal deformation 
for the supersingular elliptic curve $y^2=x^3+ax+b$ with $j$-invariant $j_k$.
The base change $t'\mapsto t'^{p-1}$ yields the family of elliptic curves
$y^2=x^3+ax+b+t'^{p-1}$, 
so its Frobenius pullback $y^2=x^3+a^px+b^p+t'^{p(p-1)}$ contains
a rational $p$-division point with osculation number one.
Making the substitution $x=t'^{-2}x'$, $y=t'^{-3}y'$, we obtain
the new Weierstrass equation $y^2=x^3+t'^{4}a^px +(b^p+t'^{p(p-1)})t'^6$,
which defines a relative cubic over $k[[t']]$. This 
cubic is already defined over the subring $k[[t]]\subset k[[t']]$,
where $t=t'^2$. According to Proposition \ref{osculation decreases},
the corresponding elliptic curve
$$
y^2=x^3+t^2a^px +(b^p+t^{p(p-1)/2})t^3
$$
over $K=k((t))$ contains a rational $p$-division point whose specialization
in the N\'eron model is nonzero.
Its reduction type is ${\rm I}_0^*$  by Proposition \ref{additive reduction}.

Write $p=1-2d$ for some integer $d$.
Making the substitution $x=t^{2d}x'$, $y=t^{3d}y'$, we obtain
the Weierstrass equation
$$
y^2=x^3+t^{2-4d}a^p x +(b^p+t^{p(p-1)/2})t^{3-6d},
$$
which is the Frobenius pullback of (\ref{generic supersingular}).
According to Proposition \ref{frobenius table}, the curve $E_K$ has reduction type
${\rm I}_0^*$.
\qed

\begin{remark}
In the proofs of Propositions
\ref{congruent 2}, \ref{congruent 3} and \ref{congruent arbitrary}
we constructed our examples from the base change 
$t'\mapsto t'^{p-1}$ from a versal deformation of 
a supersingular elliptic curve.
If we apply an $n$-fold Frobenius pullback to this base change
and carry out the constructions explained in the proofs of these
propositions, we obtain all the examples of this
section, but now with osculation number $n$.
In particular, Frobenius pullbacks from the curves constructed
in this section give all possible reduction types with arbitrary
osculation numbers.
We leave it to the reader to determine explicit Weierstrass 
equations.
\end{remark}

\section{The elliptic curve over the Igusa curve}
\mylabel{igusa neron}

Let $p>0$ be a prime number, and consider the ordinary part of the \emph{Igusa stack} $\Ig(p)^\ord$,
whose objects over a $k$-algebra $A$ are pairs $(E,z)$, where $E\ra\Spec(A)$ is
a family of ordinary elliptic curves, and $z:\Spec(A)\ra E^{(p)}$ is a section whose
fibers are  points of order $p$.
We have a commutative diagram of algebraic stacks
$$
\begin{xy}
\xymatrix{
 & \Ig(p)^\ord \ar[dl]_{\text{forget $z$}}\ar[dr]^j &\\
 \overline{M}_{1,1} \ar[rr]_j & &\PP^1
}
\end{xy}
$$
where the map on the left maps is $(E,z)\mapsto (E,0)$, which is a cyclic Galois covering of degree $p-1$,
and the map on the right is $(E,z)\mapsto j(E)$, which has degree $(p-1)/2$. Note that the horizontal map
of algebraic stacks has degree $1/2$, and that the image of the $j$-map $\Ig(p)^\ord\ra\AA^1$ is precisely the ordinary locus on the $j$-line.

This moduli problem has been first studied by Igusa \cite{Igusa 1968}.
For $p\geq 3$, the Igusa stack is representable by \cite{Katz; Mazur 1985}, Corollary 12.6.3, 
and we shall assume $p\geq 3$ in this section.
Abusing notation we  write $\Ig(p)^\ord$ for the corresponding algebraic curve,
and denote by $\Ig(p)$   its normal compactification. 
Let $U^\ord\ra\Ig(p)^\ord$ be the universal elliptic curve, and $U\ra\Ig(p)$ be
its N\'eron model.  
We now give a complete description of the universal curve around
supersingular points.
Let $F=\O_{\Ig(p),\eta}$ be the function field of the Igusa curve.

\begin{theorem}
\mylabel{main result}
Suppose $p\geq 5$. Let $x\in\Ig(p)$ be a supersingular point. 
Then a Weierstrass equation for $U_F$ over the completion at
$x\in\Ig(p)$, as well as the reduction type for 
$U_F$ and its Frobenius pullback are as in the following table:
\vspace{-0.8em}
$$
\begin{array}[t] {|*{5}{c|}}
\hline
j(x) & p & \text{\rm Weierstrass equation} & U_F & U_F^{(p)}\\
\hline
0 & \cong -1 \!\!\!\mod 3 & y^2=x^3+t^{(p-5)/6}x + t^{-1} & {\rm II}^* & {\rm II}\\
\hline
1728 & \cong -1\!\!\!\mod 4 &y^2=x^3+t^{-1}x+t^{(p-7)/4} & {\rm III}^* & {\rm III}\\
\hline
\neq 0,1728 & \text{\rm all $p$} &y^2=x^3+at^{-2p}x+(b+t^{(p-1)/2})t^{-3p} & {\rm I}^*_0 & {\rm I}^*_0\\
\hline
\end{array}
$$
Here  $t\in{\O}_{\Ig(p),x}^\wedge$ is a suitable uniformizer, and the scalars $a,b\in k$ 
in the last row are so that the elliptic curve
$y^2=x^3+ax+b$ has $j$-invariant $j(x)$. Moreover, the rational  $p$-division points in
$U_F^{(p)}$ have nonzero specialization in the N\'eron model.
\end{theorem}

\proof
Note   that the entries for the Frobenius pullback $U_F^{(p)}$ are determined by  those for $U_F$ and vice versa, according to Proposition \ref{frobenius table}.
We now  give a complete proof for the case $j(x)=0$, the other cases being   analogous  
and left to to the reader. 
So suppose $j(x)=0$. Obviously, this $j$-value must be supersingular,
whence $p\cong -1$ modulo $3$.
According to Proposition \ref{congruent 2}, there is
an elliptic curve $E_K$ over the function field $K=k((t))$ of $R=k[[t]]$,
with reduction type ${\rm II}^*$ and Weierstrass equation as in the table.
Moreover, the Frobenius pullback $E_K^{(p)}$ contains a rational  $p$-division point whose
specialization  in  the closed fiber of the N\'eron model is nonzero.
Let $\varphi:\Spec(K)\ra\Ig(p)$ be the corresponding classifying morphism,
such that $E_K=U\otimes_F K$ and $E^{(p)}_K=U^{(p)}_F\otimes K$. 
By the valuation criterion, the morphism extends uniquely to a morphism $\varphi:\Spec(R)\ra\Ig(p)$.
By Igusa's result \cite{Katz; Mazur 1985}, Corollary 12.6.2 the morphism $j:\Ig(p)\ra\PP^1$ is 
totally ramified over the supersingular $j$-values and hence
the point $x$ lies in the image of $\varphi$. 
We thus obtain an extension
of discrete valuation rings $R'\subset R$, say with ramification index $e\geq 1$, where $R'=\O_{\Ig(p),x}$. We claim that $e=1$.
To see this, denote by $\nu\geq 2$   the valuation of a minimal discriminant for $U_F^{(p)}$ at $x$.
Since $E_K^{(p)}$ has reduction type ${\rm II}$, Lemma \ref{ramification index} below tells us that $2=e\nu$,
and whence $e=1$.  Since N\'eron models are preserved under extensions with ramification index $e=1$,
the curve $U_F$  has reduction type ${\rm II}^*$ at $x$,
and the statement about the specialization of the rational point on $E_K^{(p)}$ follows
in a similar way.
%
\qed

\medskip
We have used the following observation:
Suppose $R'\subset R$ is an extension of discrete valuation rings
of arbitrary characteristic $p>0$, with the same residue field $k=R/\maxid_R=R'/\maxid_{R'}$,
and function fields $K'\subset K$.
Suppose $E_{K'}$ is an elliptic curve containing a rational $p$-division point.
Set $E_{K}=E_{K'}\otimes K$, and let  
$E\ra\Spec(R)$ and $E'\ra\Spec(R')$ be the N\'eron models of $E_K$ and $E_{K'}$, respectively.

\begin{lemma}
\mylabel{ramification index}
Suppose the rational $p$-division point on  $E_{K'}$ specializes into ${E'_k}^{0}$,
and that the induced point on $E_K$  specializes into
a nonzero element of $E_k^0$.
Let $\nu$ and $\nu'$ be the  valuations of minimal discriminants
for $E_K$ and $E_{K'}$, respectively, and $e\geq 1$ be the ramification index
of $R\subset R'$. Then we have $\nu=e\nu'$.
\end{lemma}

\proof
We have $\nu=e\nu'-12c$, where $c\geq 0$ is the number of cycles
needed in the Tate Algorithm before termination.
Consider the canonical homomorphism of relative group schemes $E'\otimes_{R'}R\ra E$.
If $c\geq 1$, then the connected component of the origin in $E'_k$ is mapped to the origin in $E_k$,
whence the rational  $p$-division point on $E_K$ specializes to zero, contradiction.
\qed

\medskip
It remains to determine the N\'eron model over the cusps of $\Ig(p)$, that
is, the points where the $j$-invariant has a pole
\cite{Katz; Mazur 1985}, Section 8.6.3.

\begin{theorem}
 Suppose $p\geq3$.
 Then the scheme of cusps of $\Ig(p)$ is finite \'etale of length $(p-1)/2$
 with a transitive action of the Galois group of $j:\Ig(p)\to\PP^1$.
 The N\'eron model over a cusp of $\Ig(p)$ has multiplicative reduction of type $I_1$
 and its Frobenius pullback has multiplicative reduction of type $I_p$. 
\end{theorem}

\proof
Since the $j$-invariant has negative valuation, $U_F$ has potentially
multiplicative reduction.
If $U_F$ had additive reduction then we would have reduction of type
 $I_\ell^*$ for some $\ell\geq1$.
However, this is excluded by Corollary \ref{restriction i} and
$U_F$ has already multiplicative reduction.

Let $K=k((t))$ and $q=t^p$.
Then there exists an elliptic curve $E_K$ over and a homomorphism
$K^*\to E_K(K)$ with kernel $q^\ZZ$, namely the Tate curve
(\cite{Silverman 1994}, Chapter V, \S3).
In particular, $t\in K^*$ maps to a rational $p$-division point of $E_K$
and since $\nu(j)=-\nu(q)=-p$ this elliptic curve has multiplicative
reduction of type $I_p$.

This curve is the Frobenius pullback of a curve induced from $U_F$
around a cusp $x\in\Ig(p)$.
Hence $\nu(j(x))=-1$ for this cusp $x$, which implies that
$U_F$ has multiplicative reduction of type $I_1$ at $x$.
Hence $j:\Ig(p)\to\PP^1$ is \'etale around $x$ and since
$j$ is a Galois morphism the same is true for every cusp.
In particular, the scheme of cusps is \'etale of length $(p-1)/2$
and $I_1$ is the reduction type of the N\'eron model for 
every cusp of $\Ig(p)$.
\qed

\medskip
Let $E$ be the N\'eron model of an elliptic curve $E_K$ and assume that
it has multiplicative reduction.
If there exists a rational $p$-division point on $E_K$ then it generates a group
scheme $G\subset E$ with generic fiber $(\ZZ/p\ZZ)_K$, which can
specialize to $\alpha_p$ or $\ZZ/p\ZZ$ only.
Since $E_k^0\iso\GG_m$ neither of these latter group schemes is contained
in $E_k^0$ and so the rational $p$-division point specializes non-trivially 
into $\Phi_k$.

\begin{remark}
 Note that Ulmer \cite{Ulmer 1990a}, Section 2 
 gave Weierstrass equations with 
 coefficients in $F$ for the universal curve $U_F$,
 which rely on relations between Eisenstein series and are of somewhat
 implicit nature. 
 There it is also shown that the universal family over $\Ig(p)$ 
 descends to the $j$-line if and only if $p\equiv -1$ modulo $4$.
 In this case, the reduction type of the N\'eron model of the 
 descended family has been determined in loc. cit., Section 6.
 It is interesting to note that the universal family
 over $\Ig(p)$
 has good reduction over $j=0$ if $p\equiv 1$ modulo $3$,
 whereas the descended family acquires additive reduction.
\end{remark}

\section{Elliptic curves with $\delta=1$}
\mylabel{delta one}

The Igusa stack in characteristic two and three has entirely new features because
wild ramification shows up.
In this section we briefly recall some relevant facts from
ramification theory  
and analyze the Galois representation on torsion points
attached to elliptic curves whose wild part of the conductor
is nontrivial yet as small as possible, namely 
$\delta=1$. These results will be applied to universal families
over Igusa curves in the next sections.
For more background on ramification groups, we refer to the 
monographs \cite{Serre 1979} and \cite{Fesenko; Vostokov 1993}
and the survey article \cite{Sautoy; Fesenko 2000}.

Suppose $k$ is an algebraically closed field
of characteristic $p>0$ and  set  $R=k[[t]]$ and $K=k((t))$.
Let $K\subset L$ be a finite Galois extension, with Galois group $G=\Gal(L/K)$,
and  $R_L\subset L$ be the integral closure of $R$.
The \emph{higher ramification groups}
$$
G_0\supset G_1\supset G_2\supset\ldots
$$
are defined as the subgroups $G_i\subset G$ that act trivially
on the $i$-th infinitesimal neighborhood $\Spec(R_L/\maxid_L^{i+1})$.
Then $G=G_0$, the $G_i\subset G$ are normal,   $G_1\subset G$ is the Sylow $p$-subgroup,
and $G/G_1$ is cyclic of order prime to $p$.
Using the existence of Sylow subgroups in $G$ for the prime divisors of   $[G:G_1]$,
we infer that $G\simeq G_1\rtimes C_m$ for some integer prime to $p$. 
Here and throughout, $C_m$ denotes a cyclic group of order $m$.

Choose a prime $l\neq p$.
The \emph{Swan representation} $P$ attached to the Galois group $G$
is the projective $\ZZ_l[G]$-module
whose character is given by $b(\sigma)=-\max\left\{i\mid \sigma\in G_i\right\}$, $\sigma\neq e$,
and $\sum_{\sigma\in G} b(\sigma)=0$.
If $V$ is a $\FF_l[G]$-module, one  defines an integer invariant
$\delta(V)=\dim_{\FF_l}\Hom_G(P,V)$, which does not depend on the choice of $K\subset L$.
It also satisfies the formula
\begin{equation}
\label{hilbert formula}
\delta(V)=\sum_{i\geq 1}\frac{1}{[G:G_i]}\dim_{\FF_l} V/V^{G_i}.
\end{equation}

Now let $E_K$ be an elliptic curve and $E$ be its N\'eron model.
Choose a Galois extension $K\subset L$ 
so that the $\FF_l$-vector space $E[l](L)$ becomes $2$-dimensional.
The invariant $\delta=\delta(E[l](L))$ is called the \emph{wild part of the conductor}.
It does not depend on $l$.
If $E_K$ has additive reduction, Ogg's formula tells us
$\nu(\Delta)=2+\delta+(m-1)$,
where $m$ denotes the number of irreducible components
in the closed fiber of the minimal model.

By construction, the Galois group $G$ comes along with a   representation on the vector space
$E[l](L)$,
which we regard as a homomorphism $G\ra\GL(2,\FF_l)$. The linear $G$-action 
respects the Weil pairing $\Phi:E[l]\times E[l]\ra\mu_l(K)$ in the sense that $\Phi(a^\sigma,b^\sigma)=\Phi(a,b)^\sigma$. Since $k$ is algebraically closed,
the $G$-action on $\mu(K)$ is trivial, such that we have a factorization $G\ra\SL(2,\FF_l)$.
We remark in passing  that this factorization holds for $l=2$ without any assumption on $k$.
Note that saying that $G\ra\SL(2,\FF_l)$ is surjective means that  the scheme of nonzero $l$-torsion $E_K[l]-0$ is connected,
and stays so under base extension as long as possible.

It is often convenient to replace $G$ by its image in $\SL(2,\FF_l)$, such that $K\subset L$ becomes
the smallest Galois extension so that $E[l](L)$ is $2$-dimensional. But it is 
useful to work with the general situation   when it comes to base change:

\begin{lemma}
\mylabel{relatively prime}
Let $K\subset K'$ be a field extension of degree $d$  prime to $p$.
Then the wild part of the conductor for the induced elliptic curve $E_{K'}$ 
is $\delta'=d\delta$.
\end{lemma}

\proof
Enlarging $L$, we may assume that $K\subset K'\subset L$, and set $G'=\Gal(L/K')$.
Obviously $G_i'=G'\cap G_i$ are the ramification groups for $K'\subset L$.
By Puiseux's Theorem, $K\subset K'$ is cyclic, hence corresponds to
a surjection $G\ra C_d$, whose kernel equals $G'$, and contains $G_1$ because $d$ is prime to $p$.
We conclude $G'_i=G_i$ for $i\geq 1$, and the statement follows from Formula (\ref{hilbert formula}).
\qed

\begin{lemma}
 Let $K\subset K'$ be a finite and purely inseparable field extension.
 Then the wild part of the conductor for the induced elliptic curve $E_{K'}$
 is $\delta'=\delta$.
\end{lemma}

\proof
Choose a $R_K$-algebra generator $x\in R_L$.
Given $\sigma\in G$ we have $\sigma\in G_i$ if and only if $\nu_L(\sigma(x)-x)\geq i+1$
by \cite{Serre 1979}, Chapter IV, \S1, Lemma 1.
Recall that $K=k((t))$ so that $K'=K^{1/p^n}$ for some $n\geq1$.
It suffices to treat the case $n=1$.
Then $L'=L\otimes_K K'$ and a straightforward argument shows that 
$x^{1/p}$ is a $R_{K'}$-algebra generator of $R_{L'}$.
Now
$$
\nu_{L'}(\sigma(x^{1/p})-x^{1/p})=\frac{1}{p}\nu_{L'}(\sigma(x)-x)=\frac{p}{p}\nu_L(\sigma(x)-x).
$$
Using this equation we conclude that the higher ramification groups and their indices 
for $L/K$ and $L'/K'$ coincide.
The statement now follows from Formula (\ref{hilbert formula}).
\qed

\medskip
The group $\SL(2,\FF_2)=\GL(2,\FF_2)$ has order $6$, consequently $\delta =0$ for characteristic $p\geq 5$.
For the rest of the section, we work in characteristic two and three
and examine the Galois representation $G\ra\SL(2,\FF_l)$ for
elliptic curves with $\delta=1$.

We start with the case $p=3$ and choose $l=2$.
Note that the action on $\PP^1(\FF_2)$ gives a bijection $\GL(2,\FF_2)\ra S_3$,
and these groups are  isomorphic to the nontrivial semidirect product $C_3\rtimes C_2$.
Let $E_K$ be an elliptic curve, and choose a Galois extension $K\subset L$ so
that $E[2](L)$ becomes a 2-dimensional $\FF_2$-vector space.
Let $\delta$ be the wild part of the conductor for $E_K$ and $g$ be the order
of the Galois group $G=\Gal(L/K)$.

\begin{proposition}
\mylabel{delta 3}
If $\delta=1$, then the homomorphism $G\ra\SL(2,\FF_2)$ is surjective.
If moreover $9\nmid g$, then $G$ is isomorphic to the
nontrivial semidirect product $G=C_3\rtimes C_{g/3}$, and  the ramification groups  are 
$G_0=G$ and  $G_1=\ldots=G_m=C_3$ and $G_{m+1}=1$  with $m=g/6$.
\end{proposition}

\proof
Suppose  that the homomorphism in question is not surjective. Replacing $G$ by its image, we may
assume that $G\subsetneq \SL(2,\FF_2)$ is a subgroup.
If $g$ is prime to $p=3$ then $\delta=0$, contradiction.
Suppose $g=3$. Since each matrix in $\GL(2,\FF_2)$ of order three
is conjugate to $(\begin{smallmatrix} 0&1\\ 1& 1\end{smallmatrix})$, we have
$E[2](L)^{G_i}=0$ for every nontrivial $G_i$, and Formula (\ref{hilbert formula}) yields
$1=\delta\geq 1/1\cdot 2$, again a contradiction. 
We conclude that $G\ra\SL(2,\FF_2)$ is surjective.

Now suppose that $9\nmid g$, so that $G_1=C_3$ is the unique Sylow $3$-subgroup.
Then $G=C_3\rtimes C_{g/3}$ is a semidirect product, which must be the nontrivial one because
$G\ra\GL(2,\FF_2)$ is surjective.
It remains to determine the orders of the ramification groups:
We have $g_0=g$ and $g_1=g_2=\ldots =g_m=3$ and $g_{m+1}=1$ for some $m\geq 1$.
Formula (\ref{hilbert formula}) yields $\delta= m\cdot\frac{1}{g/3}\cdot 2$, and the result follows.
\qed

\medskip 
Now choose a Weierstrass equation $y^2+a_1xy+a_3y=x^3+a_2x^2+a_4x+a_6$ for the elliptic curve $E_K$,
and let $K\subset K'$ be the field extension obtained by adjoining a root of
the cubic $x^3+(a_1^2+a_2)x^2+(a_4-a_1a_3)x+(a_3^2+a_6)$.

\begin{corollary}
\mylabel{change 3}
Notation as above.
If the elliptic curve $E_K$ has $\delta=1$, then $K\subset K'$ is a non-Galois extension of degree three,
the induced elliptic curve $E_{K'}$ has $\delta'=0$, and the $\FF_2$-vector space $E[2](K')$
is $1$-dimensional.
\end{corollary}

\proof
We may chose $K\subset L$ so that its Galois group is $G=\GL(2,\FF_2)$, by Proposition \ref{delta 3}.
By the inversion formula (\cite{Silverman 1986}, III.2.3.) 
the scheme of nonzero $2$-torsion on $E_K$ is given by
$y=a_1x+a_3$ together with the Weierstrass equation, whence by the cubic in question.
The Galois correspondence implies that  $K'\subset L$, and that $K\subset K'$ is non-Galois of degree three.
Using that $K'\subset L$ has degree two, we infer that $\delta'=0$ and that $E[2](K')$ is $1$-dimensional.
\qed

\medskip
For the rest of this section we work in characteristic $p=2$ and choose $l=3$.
This case is   more challenging.
To start with, let us briefly recall some well-known facts on the group $\SL(2,\FF_3)$.
We have a commutative diagram
$$
\begin{CD}
1 @>>> C_2 @>>>  \GL(2,\FF_3) @>>>  S_4 @>>>  1\\
& & @A=AA   @AAA   @AAA  \\
1 @>>> C_2 @>>> \SL(2,\FF_3)@>>> A_4@>>> 1\\
&& @A=AA   @AAA   @AAA  \\
1@>>> C_2 @>>>  Q@>>> V@>>> 1,
\end{CD}
$$
where $\GL(2,\FF_3)\ra S_4$ is  given by the action on $\PP^1(\FF_3)$.
The group $V\subset A_4$ is the Klein four group, and
$Q=\left\{\pm 1,\pm i, \pm j,\pm k\right\}$ is the quaternion group.
Its six elements of order four correspond to the traceless matrices in $\SL(2,\FF_3)$.

Since $A_4=V\rtimes C_3$ it follows that $V\subset A_4$ is the commutator subgroup,
and $S_4/V=S_3$ implies that $A_4\subset S_4$ is the commutator subgroup.
Using that $-1\in Q$ is a commutator, we infer that
$1\subset C_2\subset Q\subset\SL(2,\FF_3)\subset\GL(2,\FF_3)$
is the derived series. Let us depict the lattice of subgroups in $\SL(2,\FF_3)$:
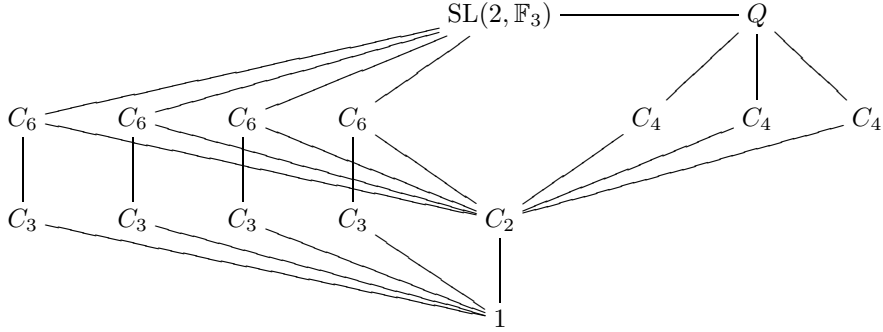
$$
\begin{xy}
\xymatrix{
&&&& \SL(2,\FF_3) \ar@{-}[dllll]\ar@{-}[dlll]\ar@{-}[dll]\ar@{-}[dl]\ar@{-}[rr] && Q\ar@{-}[dl]\ar@{-}[d]\ar@{-}[dr]\\
C_6\ar@{-}[d]\ar@{-}[drrrr] & C_6\ar@{-}[d]\ar@{-}[drrr] & C_6\ar@{-}[d]\ar@{-}[drr] & C_6\ar@{-}[d]\ar@{-}[dr] && 
C_4\ar@{-}[dl] & C_4\ar@{-}[dll] & C_4\ar@{-}[dlll]\\
C_3\ar@{-}[drrrr]&C_3\ar@{-}[drrr]&C_3\ar@{-}[drr]&C_3\ar@{-}[dr]&  C_2\ar@{-}[d]\\
&&&& 1
}
\end{xy}
$$
The   normal subgroups in $\SL(2,\FF_3)$ are precisely
$1\subset C_2\subset Q\subset\SL(2,\FF_3)$. Note that the $C_6\subset\SL(2,\FF_3)$
are the four Borel subgroups, that is, conjugate to the group 
$(\begin{smallmatrix} *&*\\0&*\end{smallmatrix})$.

Now suppose $E_K$ is an elliptic curve. Choose a Galois extension $K\subset L$
so that $E_K[3](L)$ becomes $2$-dimensional, and let $G=\Gal(L/K)\ra\SL(2,\FF_3)$
be the associated representation on $E_K[3](L)$.
The lattice of subgroups in $\SL(2,\FF_3)$
is related to a Weierstrass equation $y^2+a_1xy+a_3y = x^3+a_2x^2+a_4x+a_6$
as follows: The  subgroups $C_6\subset\SL(2,\FF_3)$ correspond to the field
extensions of $K$ obtained by adding one of the four roots of the quartic
$x^4+b_2x^3+b_4x^2+b_6x+b_8$, which defines the $x$-coordinate for one of the four lines
of points of order three 
(compare with the duplication formula \cite{Silverman 1986}, III.3.2.).
The normal subgroup $Q\subset \SL(2,\FF_3)$ corresponds to the Galois extension
of $K$ obtained by splitting the \emph{resolvent cubic} $x^3+b_4x^2+b_2b_6x+b_2^2b_8+b_6^2$
(compare \cite{Morandi 1996}, Section III.13).
Note that the quartic has Galois group contained in $V=Q/C_2$ after splitting the 
resolvent cubic.
The normal subgroup $C_2$ corresponds to splitting the
quartic. The inclusion $1\subset C_2$ finally corresponds to adding the $y$-coordinates
of the $3$-torsion points.

\begin{proposition}
\mylabel{delta 2}
Suppose the elliptic curve $E_K$ has $\delta=1$. Then  
$G\ra\SL(2,\FF_3)$ is surjective.
If moreover $16\nmid g$, then $G$ is isomorphic to the nontrivial semidirect
product $Q\rtimes C_{g/8}$, and the
ramification groups are $G_0=G$ and $G_1=...=G_s=Q$, 
$G_{s+1}=\ldots=G_{3s}=C_2$ with $s=g/24$
and $G_m=1$ for $m\geq 1+g/8$.
\end{proposition}

\proof
For the first statement, we may assume $G\subset\SL(2,\FF_3)$,
and our task is to show that $G$ has order $g=24$.
Note first that $2\mid g$, because otherwise $\delta=0$.
We observe that $(\begin{smallmatrix} -1&0\\0&-1\end{smallmatrix})\in G_i$  
and whence $E[3](L)^{G_i}=0$ for every nontrivial $G_i$.
Second, we have $3\mid g$. Otherwise $G$ would be a $2$-group,
and Formula (\ref{hilbert formula}) gives $\delta\geq 1/1\cdot 2$, contradiction.
Third, we have $g\neq 6$. Otherwise the orders of the ramification groups
are $g_1=\ldots=g_m=2$ and $g_{m+1}=1$ for some $m\geq 1$,
and thus $\delta=m\cdot 1/3\cdot 2$, contradiction.
Fourth and last, we have $g\neq 12$, because $\SL(2,\FF_3)$ contains
no subgroup of order $12$; such a group would be normal, and whence
define a splitting for $\SL(2,\FF_3)\ra A_4$, which is absurd.
Thus, $G$ is equal to $\SL(2,\FF_3)$.

Let us determine the higher ramification groups for the special case
$G=\SL(2,\FF_3)$.
Since $C_4$ is not normal in $Q$, it cannot be among the ramification groups.
We obtain $G_1=Q$, $G_2=G_3=C_2$ and $G_m=1$ for $m\geq4$ 
as in the proof for Proposition \ref{delta 3}.

Now let $G$ be arbitrary with $16\nmid g$ and denote by $G'$ the kernel of the 
surjective homomorphism $G\ra\SL(2,\FF_3)$.
The lower filtration on ramification groups does not behave well with respect to
passing to quotients. But the so-called \emph{upper filtration} $G^x=G_{\psi(x)}$
has precisely the property $(G/G')^x=G^x G'/G'$, see
\cite{Serre 1979} Section IV.\S3, Proposition 14.
Here $\psi:[0,\infty]\ra[0,\infty]$ is a convex piecewise linear homeomorphism
called the \emph{Hasse--Herbrand} function.
Its inverse $\varphi:[0,\infty]\ra[0,\infty]$ can be defined in terms of the indices of the ramification groups by
$$
\varphi(x)= 1/[G_0:G_1] + \ldots + 1/[G_0:G_n] + (x-n)/[G_0:G_{n+1}],\quad
n\leq x\leq n+1. 
$$
Knowing the ramification groups of $G/G'=\SL(2,\FF_3)$ already, 
and using the Hasse--Herbrand function, one computes the 
desired ramification groups for $G$ as in the statement.
\qed

\medskip
We now can compute the behavior of the wild part of the conductor under base changes contained in $L$:

\begin{proposition}
\mylabel{change 2}
Suppose the wild part of the conductor for  $E_K$ is $\delta=1$.
Let $G'\subset \SL(2,\FF_3)$ be  a subgroup,
and $K'\subset L$ be the   fixed field of the preimage of $G'$ in $G$. Then the wild
part of the conductor for the induced elliptic curve $E_{K'}$ is
given by the following table:
$$
\begin{array}{|*8{c|}}
\hline
G' & \SL(2,\FF_3) & C_6 & Q & C_4 & C_3 & C_2 & 1\\
\hline
\delta' & 1 & 2 &  3 & 4 &  0 &  6 & 0\\
\hline
\end{array}
$$
\end{proposition}

\proof
We may assume $G=\SL(2,\FF_3)$.
The ramification groups for $G'$ are $G'_i=G'\cap G_i$,
and the statement follows by an elementary computation from Proposition \ref{delta 2}
together with the formula (\ref{hilbert formula}).
Consider, for example, the case $G'=C_4$.
Then we have $G'_0=C_4$ and  $G'_1=C_4$, $G'_2=G'_3=C_2$ and $G'_4=0$,
consequently $\delta'= 1/1\cdot 2+ 1/2\cdot 2+ 1/2\cdot 2= 4$.
\qed

\medskip
Now let $G'\subset G\subset \SL(2,\FF_3)$ be two subgroups
so that $G'\subset G$ has index two.
Then the extension of fixed fields  $F\subset F'$ is cyclic of degree two.
Let $R_F\subset R_{F'}$ be the corresponding extension of discrete valuation rings.
To control Weierstrass equations, it will later be useful to express the uniformizer of $R_{F}$ in terms
of the uniformizer of $R_{F'}$. Luckily, the situation is as simple as possible:

\begin{proposition}
\mylabel{quadratic uniformizer}
Suppose that $E_K$ has $\delta=1$, and let $u\in R_F$ be a uniformizer.
Then there is a uniformizer $s\in R_{F'}$ and a nonzero scalar $\lambda\in k$ with $u=s^2/(\lambda-s)$.
\end{proposition}

\proof
The isomorphism class of $F\subset F'$ corresponds to an element from the cohomology group
 $H^1(F,\ZZ/2\ZZ)=F/\wp(F)$, $\wp(g)=g^2-g$, which we identify  with the group
of odd polynomials $f$ in $u^{-1}$. Let $f=\lambda_{1-2n}u^{1-2n}+\ldots+\lambda_{-1}u^{-1}$
be the odd polynomial for $F\subset F'$, with $n\geq 1$ and $\lambda_{1-2n}\neq 0$.
Then $R_{F'}=R_F[s]/(s^2-u^ns -u^{2n}f)$,
and the Galois involution is $s\mapsto s+u^n$. 
Since $u=s^2+O(3)$, the ramification groups for $G/G'$ are
$$
(G/G')_0=\ldots =(G/G')_{2n-1}=C_2\quadand (G/G')_{2n}=1.
$$
We have to show $n=1$, or, equivalently  $(G/G')_2=1$.
As explained in the proof of Proposition \ref{delta 2}, we have to
pass to the upper filtration and determine the Hasse--Herbrand function.
In our situation, the Hasse--Herbrand function for $G/G'$ has a unique
break point at $x=2n-1$.

We now make the computation in the case $G=C_4$, $G'=C_2$, and leave the
other cases to the reader.
By Proposition \ref{delta 2}, the ramification groups for $G$ are
$G_0=C_4$, $G_1=G_2=G_3=C_2$ and $G_4=1$. We infer
that the Hasse--Herbrand function for $G/G'$ has its break point
at $x=1$, and consequently $n=1$.
\qed

\section{The Igusa curve in characteristic three}
\mylabel{igusa 3}

In this section, we shall analyze the the  N\'eron model $U\ra\Ig(3)$ and  the resulting 
possibilities for rational points of order three with nonzero specialization in the N\'eron model.
First note that the $j$-map $j:\Ig(3)\ra\PP^1$,
having degree $1=(3-1)/2$, is an isomorphism. Second note that there is only one supersingular point $x\in\Ig(3)$, which has $j$-value $j(x)=0$.
So we may identify $\Ig(3)^\ord=\AA^1-\left\{0\right\}=\Spec(A)$, where $A=k[t^{\pm1}]$ is the ring of Laurent polynomials,
and $t$ is a uniformizer at $j=0$.
Since $\Pic(A)=0$, the  universal elliptic curve $U^\ord$ 
must admit a global Weierstrass equation over $A$:

\begin{proposition}
\mylabel{universal 3}
The universal elliptic curve $U^\ord$ has as global Weierstrass equation
$y^2+xy=x^3-1/t$ over $A$.
At the supersingular point $x\in\Ig(3)$, the 
reduction type   is ${\rm II}^*$, the valuation of a minimal discriminant
is $\nu(\Delta)=11$, and the wild part of the conductor is $\delta=1$.
\end{proposition}

\proof
The given Weierstrass equation has $j$-invariant $j=t$ and discriminant $1/t$.
Whence it differs from the universal elliptic curve by a quadratic twist.
Using inversion and duplication formula, we see that
the Frobenius pullback $y^2+txy=x^3-t^3$ admits a rational point
of order three, namely $x=t$, $y=0$.
Since a nontrivial quadratic twist   destroys this rational point,
we conclude that $U^\ord$ is actually given by this Weierstrass equation.
The remaining statements follow from the Tate Algorithm and  Ogg's Formula.
\qed

\begin{remark}
 The universal curve over $\Ig(3)$ has already been determined in 
 \cite{Ulmer 1990a}, Proposition 2.3.
 For the sake of completeness we decided to include a proof in our setup.
\end{remark}

\medskip
Now set $R=k[[t]]$ and consider the induced elliptic curve $U_K$ over $K=k((t))$.
Since $\delta=1$, our results form the preceding section apply.
The goal now is to construct elliptic curves  so that the  rational points of order three on the Frobenius pullbacks
have nonzero specialization in the N\'eron model.
We   can compute the behavior of $U_{K'}=U_K\otimes K'$ for various base changes $K\subset K'$   of successive degrees $d=2,3,5$.
using Lemma \ref{relatively prime} and Proposition \ref{change 2}.
Our findings are best summarized in a family tree:
$$
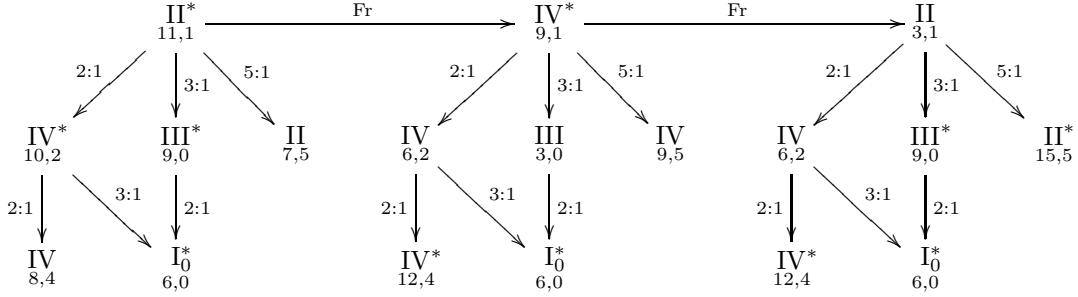
\hspace{-2em}
\begin{xy}
\xymatrix{
&
\underset{11,1}{\p{\rm II}^*} \ar[dl]_{2:1} \ar[d]^{3:1} \ar[dr]^{5:1} \ar[rrr]^{\Fr} 
&&&
\underset{9,1}{\p{\rm IV}^*}  \ar[dl]_{2:1} \ar[d]^{3:1} \ar[dr]^{5:1} \ar[rrr]^{\Fr} 
&&&
\underset{3,1}{{\rm II}} \ar[dl]_{2:1} \ar[d]^{3:1} \ar[dr]^{5:1}
\\
\underset{10,2}{\p{\rm IV}^*}\ar[dr]^{3:1} \ar[d]_{2:1}
&
\underset{9,0}{\p{\rm III}^*}\ar[d]^{2:1}
&
\underset{7,5}{{\rm II}}
&
\underset{6,2}{{\rm IV}}\ar[dr]^{3:1} \ar[d]_{2:1}
&
\underset{3,0}{{\rm III}}\ar[d]^{2:1}
&
\underset{9,5}{{\rm IV}}
&
\underset{6,2}{{\rm IV}}\ar[dr]^{3:1} \ar[d]_{2:1}
&
\underset{9,0}{\p{\rm III}^*}\ar[d]^{2:1}
&
\underset{15,5}{\p{\rm II}^*}
\\
\underset{8,4}{{\rm IV}}
&
\underset{6,0}{\p{\rm I}_0^*}
&&
\underset{12,4}{\p{\rm IV}^*}
&
\underset{6,0}{\p{\rm I}_0^*}
&&
\underset{12,4}{\p{\rm IV}^*}
&
\underset{6,0}{\p{\rm I}_0^*}
}
\end{xy}
$$
Here the two numbers below the Kodaira symbols  denote the valuation of a minimal discriminant $\nu(\Delta)$
and the wild part of the conductor $\delta$, and the $3:1$ extensions are obtained by adjoining the root
of the cubic $x^3+t^2x^2-t^5$ (compare Corollary \ref{change 3}).
Note that $s=t^2/x$, which satisfies the integral equation 
$s^3-ts-t$ or equivalently $t=s^3/(s+1)$, 
is a uniformizer for the corresponding discrete valuation ring.

It turns out that the rational $3$-division points occurring
in the Frobenius pullbacks have nonzero specialization in the N\'eron models.
To verify this, we  compute the minimal Weierstrass equations for the $U_{K'}=U_K\otimes K'$,
using the substitutions $t=t'^2$, $t=t'^5$ or  $t=t'^3/(1+t')$:
$$
\begin{array}[t]{|l |*{3}{c|}l|}
\hline
\text{minimal Weierstrass equation} & \nu(\Delta) & \delta & \text{type} & j \\
\hline
y^2 + txy = x^3 - t^5                 &  11 & 1 & \p{\rm II}^* & t\\
y^2 + txy = x^3 - t^4                 &  10 & 2 & \p{\rm IV}^* & t^2\\
y^2 + txy = x^3 - t^2                 &   8 & 4 & {\rm IV}   & t^4\\
y^2 + txy = x^3 - t                   &   7 & 5 & {\rm II}   & t^5 \\
y^2 + txy = x^3 - t^3(1+t)            &   9 & 0 & \p{\rm III}^*& t^3/(1+t)\\
y^2 + txy = x^3 - t^2x                &   6 & 0 & \p{\rm I}_0^*& t^6/(1+t^2)\\
\hline
\hline
y^2 + txy + t^2y = x^3              &  9 & 1 & \p{\rm IV}^* & t^3\\
y^2 + txy + ty = x^3                &  6 & 2 & {\rm IV}   & t^6\\
y^2 + t^2xy + t^2y = x^3            & 12 & 4 & \p{\rm IV}^* & t^{12}\\
y^2 + t^2xy + ty = x^3              &  9 & 5 & {\rm IV}   & t^{15}\\
y^2 + txy + (1+t)y = x^3            &  3 & 0 & {\rm III}  & t^9/(1+t)^3\\
y^2 + t^2xy + (1+t^2)y = x^3        &  6 & 0 & \p{\rm I}_0^*& t^{18}/(1+t^2)^3\\
\hline
\hline
y^2 + txy + y = x^3               &  3 & 1 & {\rm II} & t^9\\
y^2 + t^2xy + y = x^3             &  6 & 2 & {\rm IV} & t^{18}\\
y^2 + t^4xy + y = x^3             & 12 & 4 & \p{\rm IV}^* & t^{36}\\
y^2 + t^5xy + y = x^3             & 15 & 5 & \p{\rm II}^*& t^{45}\\
y^2 + t^3xy + (1+t)^3y = x^3      &  9 & 0 & \p{\rm III}^*& t^{27}/(1+t^9)\\
y^2+t^5xy=x^3-t^2x                &  6 & 0 & \p{\rm I}_0^*& t^{54}/(1+t^{18})\\
\hline
\end{array}
$$
For such families, it is easy to determine the specialization behavior of points of order three:

\begin{proposition}
Let $E_A$ be an elliptic curve over an arbitrary ring $A$ of characteristic three.
Then $E_A$ admits a section
whose fibers are rational $3$-division points if and only if it admits
a global Weierstrass equation  of the form $y^2+a_1xy+a_3y=x^3$
for some units $a_1,a_3\in A$. 
One such section of order three is then given by $x=y=0$.
\end{proposition}

\proof
The condition is necessary, because the Frobenius pullback of $U^\ord$ admits
the Weierstrass equation
$y^2+txy+t^2y=x^3$. The sufficiency follows from the duplication and the inversion formula.
\qed

\medskip
As an immediate consequence:

\begin{corollary}
\mylabel{criterion nonzero}
Suppose  $y^2+a_1xy+a_3y=x^3$ is a minimal Weierstrass equation with nonzero $a_1,a_3\in R$.
Then the rational $3$-division points on $E_K$
have nonzero specialization in the closed fiber of the N\'eron model.
Moreover, we have $a_3\in\maxid_R$ if and only if the
rational $3$-division points have nonzero class in $\Phi_k$.
\end{corollary}

\proof
Clearly, the point $z=(0,0)\in E_K$  of order three does not specialize to infinity in the Weierstrass model,
whence has nonzero specialization into the N\'eron model. Moreover, $z$ has nonzero class in $\Phi_k$
if and only if it specializes into the singularity of the Weierstrass model.
The latter is given by  $x=0$, $y=-a_3$, and the result follows.
\qed

\medskip
Examining our table above,
we obtain the following result:

\begin{theorem}
\mylabel{nonzero 3}
For  the Kodaira symbols ${\rm II}$, ${\rm II}^*$, ${\rm III}$, ${\rm III}^*$, ${\rm IV}$, ${\rm IV}^*$, ${\rm I}_0^*$, 
there is an elliptic curve $E_K$ containing a rational $3$-division point
with nonzero specialization in $E_k$ and the given reduction type.
For ${\rm IV}$ and ${\rm IV}^*$, there are such examples with nonzero specialization in $\Phi_k$, and
examples with zero  specialization in $\Phi_k$.
\end{theorem}

\bigskip
We close this section by discussing the elliptic curves 
$$
E_{n,K}: y^2+t^{2^n}xy+t^{2^{n+1}}y=x^3, \quad n\geq 0
$$
which contain a rational $3$-divsion point.  
They are obtained from the Frobenius pullback of the universal elliptic curve by the base change
of degree $2^n$. Let $\nu(\Delta)$  be the valuation of a minimal discriminant
for $E_{n,K}$

\begin{proposition}
\mylabel{family 3}
If $n$ is odd, then   $\nu(\Delta)=2^n+4$
and the reduction type of $E_{n,K}$ is ${\rm IV}$. If $n$ is even, then $\nu(\Delta)=2^n+8$,
and the reduction type is ${\rm IV}^*$. In any case, $\delta =2^n$,
and the rational $3$-division points have nonzero specialization in $\Phi_k$.
\end{proposition}

\proof
We have $\delta=2^n$ by Lemma \ref{relatively prime}.
Suppose $n$ is odd. Then $3\mid 2^{n+1}-1$, and  $E_K$ is
given by the integral Weierstrass equation $y^2+t^{2^n-(2^{n+1}-1)/3}xy + ty=x^3$, whose discriminant has
valuation $2^n+4$. Ogg's Formula implies that the latter Weierstrass equation is minimal,
and that the reduction type is ${\rm IV}$. Corollary \ref{criterion nonzero} shows that
the rational $3$-division points have nonzero class in $\Phi_k$.
The argument for $n$ even is similar.
\qed

\medskip
We see that the property of having a rational $3$-division point with nonzero class in $\Phi_k$
can be preserved under base changes of arbitrarily large degree.

\section{The Igusa stack in characteristic two}
\mylabel{igusa stack}

In this section, $A=k[t^{\pm 1}]$ denotes the ring of Laurent polynomials
over an algebraically closed field $k$ of characteristic $p=2$.
We shall analyze the Igusa stack $\Ig(2)\ra\Spec(A)$,
and in particular the reduction types of tautological families.
Here a  family of elliptic curves $E_A$ over $A$
is called a \emph{tautological family} if $j(E_A)=t$.
For example, the Weierstrass equation
\begin{equation}
\label{tautological family}
y^2+xy=x^3+a_2x^2+t^{-1}
\end{equation}
has $j$-invariant $j=t$ and discriminant $\Delta=1/t$,
and hence yields a tautological family. Note that this is independent of
the coefficient $a_2\in A$. 

When we regard a tautological family as an object
in the Igusa stack, we also call it a \emph{tautological object}.
The existence of tautological objects shows that the $\left\{\pm 1\right\}$-gerbe $\Ig(2)\ra\Spec(A)$
is trivial, that is, isomorphic to the classifying stack $B(\ZZ/2\ZZ)$.
In some sense, tautological objects are the best replacement, in a stack theoretical context, for the universal object.
To understand the set of tautological objects, consider the map
$$
\tau:  A\lra \Ig(2)_A,\quad
a_2\longmapsto E: y^2+xy=x^3+a_2x^2+t^{-1}
$$
from the group of polynomials into the set of tautological objects of the Igusa stack.
Let $A_\odd\subset A$ be the vector space of all odd Laurent polynomials.

\begin{proposition}
\mylabel{tautological objects}
The map $\tau$ induces a bijection between the group 
$A_\odd$ of odd Laurent polynomials and
the set of isomorphism classes of tautological objects in the stack $\Ig(2)$.
\end{proposition}

\proof
Consider the additive map  $\wp:A\ra A$, $f\mapsto f^2-f$. Using that $k$ is
algebraically closed, we easily see that the canonical projection 
$$
A_\odd\lra H^1(A,\ZZ/2\ZZ)=A/\wp(A)
$$
is bijective. To proceed, let $E_A$ be the tautological family given
by the Weierstrass equation (\ref{tautological family}).
It remains to see that given $a_2\in A_\odd$, the corresponding quadratic twist of $E_A$
is given by the   Weierstrass equation $y^2+xy=x^3+a_2x^2+t^{-1}$.
We sketch the argument:
The sign involution acts on $E_A$ via $y\mapsto y+x$, and the Galois involution acts
on $A[u]/(u^2-u-a_2)$ via $u\mapsto u+1$.
Whence $y'=y+xu$ and $x'=x$ are invariant under the diagonal action,
and indeed yield the desired Weierstrass equation.
\qed

\medskip
Our next task is to determine the reduction types at $t=0$.
Let $a_2\in A_\odd$ be an odd Laurent polynomial, and write its vanishing order at $t=0$
in the form $\nu(a_2)=-2d-1$. In other words, we have
$$
a_2 =  t^{-2d-1}f^2
$$
for some integer $d$ and some polynomial $f\in k[t]$ that  has nonzero constant term  or is the zero polynomial.
In the latter case we take it that $d=-\infty$.

\begin{proposition}
Let $E_A$ be the tautological family $y^2+xy=x^3+a_2x^2+t^{-1}$.

\begin{enumerate}
\mylabel{tautological types}
\item
If $d<0$, then the reduction type of $E_A$ at $t=0$ is
${\rm II}^*$, the valuation of a minimal discriminant is $\nu(\Delta)=11$,
and the wild part of the conductor is $\delta=1$.
\item
If $d\geq 0$, then the reduction type at $t=0$ is  
${\rm I}^*_{8d+3}$, the valuation of a minimal discriminant is $\nu(\Delta)=12d+11$,
and the wild part of the conductor is $\delta=4d+2$.
\end{enumerate}
\end{proposition}

\proof
We may replace the ring of  Laurent polynomials   by the field of formal Laurent series   $K=k((t))$.
Suppose that $d<0$, such that $a_2\in k[[t]]$. 
Since all power series vanish in
$H^1(K,\ZZ/2\ZZ)=K/\wp(K)$, $\wp(f)=f^2-f$,
we may as well assume  that $a_2=0$.
Then $y^2+txy=x^3+t^5$ is a minimal Weierstrass equation for $E_K$,
and statement (i) immediately follows from the Tate Algorithm.

Now suppose $d\geq 0$. Then $f\in k[[t]]$ is a unit; set $g=1/f$. 
Starting with the original Weierstrass equation, we make the substitution 
$x=(gt^{d+1})^{-2}x'+gt^{d}$, 
and obtain a new Weierstrass equation
$$
y^2 + (gt^{d+1})xy + (g^4t^{4d+3})y = x^3 +(t+g^3t^{3d+2})x^2 + (g^6t^{6d+4})x + (g^9t^{9d+6}).
$$
The coefficients of this Weierstrass equation satisfy the assumption of Lemma \ref{reduction 1} below,
which tells us that the last Weierstrass equation is minimal, and
that the reduction type is ${\rm I}^*_{8d+3}$. The remaining statements follow from Ogg's Formula.
\qed

\medskip
We next consider   Frobenius pullbacks of our tautological families:

\begin{proposition}
\mylabel{pullback types}
Let $E_A$ be the tautological family $y^2+xy=x^3+a_2x^2+t^{-1}$.
Then the $2$-torsion section on
the Frobenius pullback $E^{(2)}_A$ has nonzero specialization in
the component group $\Phi_k$ at $t=0$. Moreover:
\begin{enumerate}
\item
If $d<0$, then the reduction type is
${\rm III}^*$, the valuation of a minimal discriminant is $\nu(\Delta)=10$,
and the wild part of the conductor is $\delta=1$.
\item
If $d\geq 0$, then the reduction type is  
${\rm I}^*_{8d+2}$, the valuation of the minimal discriminant is $\nu(\Delta)=12d+10$,
and the wild part of the conductor is $\delta=4d+2$.
\end{enumerate}
\end{proposition}

\proof
To check (i), it suffices to treat the case $a_2=0$.
Then the Weierstrass equation $y^2+txy=x^3+t^4$  for $E^{(2)}_A$ must be minimal, because
$\nu(\Delta)=10$, and the result follows from the Tate Algorithm.

We next verify (ii). The Weierstrass equation for the Frobenius pullback is
$$
y^2 + xy = x^3 + f^4t^{2(-2d-1)}x^2 + t^{-2}.
$$
Again we  may replace the ring of Laurent polynomials   by
the field of formal Laurent  series   $K=k((t))$, and set $g=1/f$.
Applying successively the substitutions
$$
y = y' + f^2t^{-2d-1}x\quadand
x = (gt^{d+1})^{-2}x'\quadand
y = y' + g^3t^{3d+2},
$$
we simplify the coefficient of $x^2$, make the Weierstrass equation integral,
and remove the constant term, respectively. The outcome is the new Weierstrass equation
$$
y^2 + gt^{d+1}xy = x^3 + tx^2 + t^{4d+3}g^4x.
$$
Now Lemma \ref{reduction 2} below yields (ii).

It remains to prove the statement about the section of order $2$.
Using that all our minimal Weierstrass equations have $a_3=a_6=0$, we infer that 
the section of order $2$ is given by $x=y=0$, and hence specializes into
the singularity of the Weierstrass model. It follows that its specialization
into the component group $\Phi_k$ of the N\'eron model is nontrivial.
\qed

\medskip
Our final task is to compute the reduction types at $t=\infty$.
Changing notation, we write the  vanishing order of $a_2$ at $t=\infty$
in the form $\nu(a_2)=-2d-1$ for some integer $d$. In other words, we have
$a_2 =  t^{2d+1}f^2$
for some integer $d$ and some polynomial $f\in k[t^{-1}]$ 
that has nonzero constant term  or is the zero polynomial.

\begin{proposition}
\mylabel{infinity type}
Let $E_A$ be the tautological family $y^2+xy=x^3+a_2x^2+t^{-1}$.
\begin{enumerate}
\item
If $d<0$, then the reduction type of $E_A$ at $t=\infty$ is
${\rm I}_1$, the valuation of a minimal discriminant is $\nu(\Delta)=1$,
and the wild part of the conductor is $\delta=0$.
\item
If $d\geq 0$, then the reduction type at $t=\infty$ is  
${\rm I}^*_{8d+5}$, the valuation of a minimal discriminant is $\nu(\Delta)=12d+13$,
and the wild part of the conductor is $\delta=4d+2$.
\end{enumerate}
\end{proposition}

\proof
The arguments are as in the proof for Proposition \ref{tautological types}.
In case $d\geq 0$, one has to use the substitution 
$x=(gt^{-d-1})^{-2}x'+gt^{-d-1}$, 
where $g=f^{-1}$, which yields
the minimal Weierstrass equation
$$
y^2+gt^{-d-1}xy + g^4t^{-4d-4}y = x^3 + (t^{-1}+ g^3t^{-3d-3})x^2 + g^6t^{-6d-6}x + g^9t^{-9d-9},
$$
and Lemma \ref{reduction 1} yields the result; details are left to the reader.
\qed

\medskip
With the techniques presented in the proofs of Proposition \ref{infinity type} and
Proposition \ref{pullback types} we obtain:

\begin{proposition}
\mylabel{infinity pullback}
  Let $E_A$ be the tautological family $y^2+xy=x^3+a_2x^2+t^{-1}$.
  Then the $2$-torsion section on the Frobenius pullback $E_A^{(2)}$ 
  has nonzero specialization in the component group at $t=\infty$.
  Moreover
  \begin{enumerate}
  \item
  If $d<0$, then the reduction type is
  ${\rm I}_2$, the valuation of a minimal discriminant is $\nu(\Delta)=2$,
  and the wild part of the conductor is $\delta=0$.
  \item
  If $d\geq 0$, then the reduction type is  
  ${\rm I}^*_{8d+6}$, the valuation of a minimal discriminant is $\nu(\Delta)=12d+14$,
  and the wild part of the conductor is $\delta=4d+2$.
\end{enumerate}
\end{proposition}

\medskip
We summarize our results about tautological families in the following table
$$
\begin{array}[t]{|l|c|ccc|ccc|}
\hline
\text{around $t$} & d & \multicolumn{3}{c|}{E_A} & \multicolumn{3}{c|}{E_A^{(2)}} \\
&& \nu(\Delta) & \delta & \text{type} &  \nu(\Delta) &\delta& \text{type} \\
\hline
t=0      & <0   & 11     &    1 & II^*       & 10     &    1 & III^* \\
         &\geq0 & 12d+11 & 4d+2 & I_{8d+3}^* & 12d+10 & 4d+2 & I_{8d+2}^*\\
\hline
t=\infty & <0   &  1     &    0 & I_1        & 2      &    0 & I_2\\
         & \geq0& 12d+13  & 4d+2 & I_{8d+5}^*& 12d+14 & 4d+2 & I^*_{8d+6} \\
\hline
\end{array}
$$
In particular, we see that reduction of type $I_l^*$ plays a special role.
Note that in characteristic $p\geq3$ this type of reduction does not appear
at all if the Frobenius pullback has a rational $p$-division point 
by Corollary \ref{restriction i}.
We will come back to this point in Section \ref{semistable reduction 2}.

\medskip
We end the section by the following two technical observations, 
which are very useful in handling reduction type 
${\rm I}^*_l$, and have been used in preceding proofs:

\begin{lemma}
\mylabel{reduction 1}
Let  $E_K$ be an elliptic curve with
Weierstrass equation $y^2+a_1xy+a_3y=x^3+a_2x^2+a_4x+a_6$.
Set $n=\nu(a_3)$, and assume
$$
\nu(a_1)\geq 1,\quad\nu(a_2)=1,\quad n=\nu(a_3)\geq 2,\quad\nu(a_4)\geq n,
\quadand\nu(a_6)\geq 2n-1.
$$
Then the Weierstrass equation is minimal, and $E_K$ has reduction type ${\rm I}^*_{2n-3}$.
\end{lemma}

\proof
The Tate Algorithm reveals that the Weierstrass equation is minimal and  has reduction type ${\rm I}^*_l$ for some   $l\geq 0$.
To determine $l$, one first blows up the ideal $(x,y,t)$, which defines the reduced singular point,
and   then the ideal $(y/t,t)$, which defines the reduced fiber. In the chart with
coordinate $x/t,y/t^2,t$, the Weierstrass equation transforms into
$$
t(y/t^2)^2 + a_1(x/t)(y/t^2) + t^{-1}a_3(y/t^2)=
t^3(x/t)^3 + t^{-2}a_4(x/t) + t^{-3}a_6,
$$
and the reduced fiber becomes a configuration of $4$ copies of $\PP^1$ with intersection
graph $D_4$, containing one rational double point.
To proceed, one successively blows up the ideals
$$
(x/t,t),\quad (y/t^2,t),\quad (x/t^2,t),\quad (y/t^2,t),\quad \ldots\quad
(y/t^{n-1},t),\quad (x/t^{n-1},t).
$$
In other words, one blows up   reduced fibers $2n-3$ times. In all but the last blowing ups
this introduces on additional copy of $\PP^1$ into the fiber, whereas the last blowing up
adds two disjoint copies of $\PP^1$.
In each blowing ups, the coefficients $a_i'$ of the successive equations acquire   factors, which are
given by the following table:
$$
\begin{array}[t]{l |*{5}{c}}
\text{coefficients}                  & a'_1 & a'_3 & a'_2 & a'_4 & a'_6  \\
\hline
\text{blowing up of $(x/t^i,t)$}     & 1 & 1/t & t    &   1 & 1/t  \\
\text{blowing up of $(y/t^{i+1},t)$} & 1 & 1   & 1/t & 1/t & 1/t \\
\end{array}
$$
According to our assumptions on the original coefficients $a_i$, it is indeed possible
to carry out the  sequence of blowing ups.
After the last blowing up, the resulting scheme is regular, and we infer that
the reduction type of $E_K$ is ${\rm I}^*_l$, where the number of irreducible components
is  $4+(2n-4)+2=2n+2$, and therefore $l=2n-3$.
\qed

\medskip
Similar arguments yield the next result, whose proof is left to the reader:

\begin{lemma}
\mylabel{reduction 2}
Let $E_K$ be an elliptic curve over the field $K$ with Weierstrass equation
$y^2 + a_1xy + a_3y = x^3 + a_2x^2 + a_4x + a_6$. Set  $n=\nu(a_4)$, and suppose that
$$
\nu(a_1)\geq 1,\quad\nu(a_2)=1,\quad \nu(a_3)\geq n-1, \quad n=\nu(a_4)\geq 2,\quad\nu(a_6)\geq 2n-1.
$$
Then the Weierstrass equation is minimal, and $E_K$ has reduction type ${\rm I}^*_{2n-4}$.
\end{lemma}


\section{Other reduction types in characteristic two}
\mylabel{igusa 2}

We continue to work in characteristic $p=2$, with
$R=k[[t]]$ and $K=k((t))$.
In this section, we shall construct elliptic curves  
whose Frobenius pullbacks have various reduction types
and whose rational $2$-division point has nonzero specialization.

The starting point is the tautological curve $E_K$ given
by $y^2+xy=x^3+t^{-1}$, which has   minimal Weierstrass
equation $y^2+txy=x^3+t^5$, numerical invariants
$\nu(\delta)=11$, $\delta=1$, and reduction type ${\rm II}^*$.
Applying the results from Section \ref{delta one}, we
now examine various pullbacks for successive field extensions of degree $d=2,3,4,5$.
We depict our relevant findings in a family tree:
$$
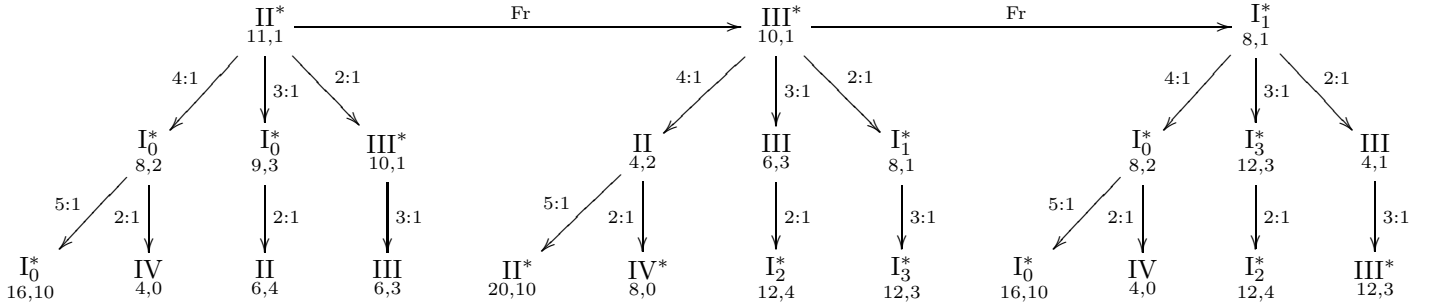
\hspace{-9em}
\begin{xy}
\xymatrix{
&&
\underset{11,1}{\phantom{{}^*}{\rm II}^*}\ar[dl]_{4:1}\ar[d]^{3:1}\ar[dr]^{2:1}\ar[rrrr]^{\Fr}
&&&&
\underset{10,1}{\phantom{{}^*}{\rm III}^*}\ar[dl]_{4:1}\ar[d]^{3:1}\ar[dr]^{2:1}\ar[rrrr]^{\Fr}
&&&&
\underset{8,1}{\phantom{{}^*}{\rm I}_1^*}\ar[dl]_{4:1}\ar[d]^{3:1}\ar[dr]^{2:1}
\\
&
\underset{8,2}{{\rm I}_0^*}\ar[dl]_{5:1}\ar[d]_{2:1}&
\underset{9,3}{\phantom{{}^*}{\rm I}_0^*}\ar[d]^{2:1}&
\underset{10,1}{{\rm III}^*}\ar[d]^{3:1}&
&
\underset{4,2}{{\rm II}}\ar[dl]_{5:1}\ar[d]_{2:1}&
\underset{6,3}{{\rm III}}\ar[d]^{2:1}&
\underset{8,1}{{\rm I}_1^*}\ar[d]^{3:1}&
&
\underset{8,2}{{\rm I}_0^*}\ar[dl]_{5:1}\ar[d]_{2:1}&
\underset{12,3}{{\rm I}_3^*}\ar[d]^{2:1}&
\underset{4,1}{{\rm III}}\ar[d]^{3:1}
\\
\underset{16,10}{{\rm I}_0^*}&
\underset{4,0}{{\rm IV}}&
\underset{6,4}{{\rm II}}&
\underset{6,3}{{\rm III}}
&
\underset{20,10}{\p{\rm II}^*}&
\underset{8,0}{\p{\rm IV}^*}&
\underset{12,4}{{\rm I}_2^*}&
\underset{12,3}{{\rm I}_3^*}
&
\underset{16,10}{{\rm I}_0^*}&
\underset{4,0}{{\rm IV}}&
\underset{12,4}{{\rm I}_2^*}&
\underset{12,3}{{\rm III}^*}
}
\end{xy}
$$
Here the two numbers below the Kodaira symbols denote the valuation of
a minimal discriminant $\nu(\Delta)$ and the wild part of the conductor $\delta$.
The quadratic extensions are given by $t=s^2/(1-s)$, and
the quartic extension is given by adjoining one root of the quartic
$x^4+t^2x^3+t^7$, which defines the $x$-coordinate for one  line 
of points of order three on $E_K$.
Note that $s=t^2/x$, which satisfies the integral equation $s^4+ts+t=0$
or equivalently $t=s^4/(1-s)$ is a uniformizer for the corresponding discrete valuation ring.

The behavior of $\delta$ follows from Proposition \ref{relatively prime},
except for the branches starting with an initial base change of degree two;
for those, $\delta$ must be computed via the Tate Algorithm.
We now can tabulate the minimal Weierstrass equations for the induced  elliptic curves
$E_{K'}=E_K\otimes K'$, using successively the substitutions
$t=t'^2/(1-t')$, $t=t'^3$, $t=t'^4/(1-t')$, and $t=t'^5$:

$$
\begin{array}[t]{|l |*{3}{c|}l|}
\hline
\text{minimal Weierstrass equation} & \nu(\Delta) & \delta & \text{type} & j \\
\hline
y^2 + txy = x^3 + t^5                 & 11 &  1 & \p{\rm II}^*  & t\\
y^2 + txy = x^3 + t^2(1+t)            &  8 &  2 & \p{\rm I}_0^* & t^4/(1+t)\\
y^2 + t^3xy + t^4y = x^3 + t^3        & 16 & 10 & \p{\rm I}_0^* & t^{20}/(1+t^5)\\
y^2 + txy + ty = x^3 + x^2 + tx + t^3 &  4 &  0 & \p{\rm IV}    & t^8/(1+t)^3(1+t+t^2)\\
y^2 + txy = x^3 + t^3                 &  9 &  3 & \p{\rm I}_0^* & t^3\\
y^2 + txy = x^3 + (1+t)^3             &  6 &  4 & {\rm II}    & t^6/(1+t^3)\\
y^2 + txy = x^3 + t^4(1+t)            & 10 &  1 & \p{\rm III}^* & t^2/(1+t)\\
y^2 + txy = x^3 + tx + t^3            &  6 &  3 & {\rm III}   & t^6/(1+t^3)\\
\hline
\hline
y^2 + txy = x^3 + t^3x                & 10 &  1 & \p{\rm III}^* & t^2\\
y^2 + txy = x^3 + (1+t)x              &  4 &  2 & {\rm II}    & t^8/(1+t)^2\\
y^2 + t^5xy = x^3 + (1+t^5)x          & 20 & 10 & \p{\rm II}^*  & t^{40}/(1+t^5)^2\\
y^2 + t^2xy = x^3 + (1+t)^2(1+t^3)x &  8 &  0 & \p{\rm IV}^*  & t^{16}/(1+t)^6(1+t+t^2)^2 \\
y^2 + txy = x^3 + tx                  &  6 &  3 & {\rm III}   & t^6\\
y^2 + t^2xy = x^3 + t^2(1+t)^3x       & 12 &  4 & \p{\rm I}^*_2 & t^{12}/(1+t)^6\\
y^2 + txy = x^3 + t^2x                &  8 &  1 & \p{\rm I}_1^*   & t^4/(1+t^2)\\
y^2 + t^2xy = x^3 + t^2x              & 12 &  8 & \p{\rm I}_3^* & t^{12}/(1+t^6)\\
\hline
\hline
y^2 + txy = x^3 + t^2x                &  8 &  1 & \p{\rm I}_1^* & t^4\\
y^2 + t^2xy = x^3 + (1+t^2)x          &  8 &  2 & \p{\rm I}_0^* & t^{16}/(1+t^4)\\
y^2 + t^8xy + t^4y = x^3 + t^6x^2+t^7x + t^3 & 16 & 10 & \p{\rm I}_0^* & t^{80}/(1+t^{20})\\
y^2 + t^3xy + ty = x^3 + t^2x^2 + t^6x+t^4   &  4 &  0 & {\rm IV} & t^{32}/(1+t)^{12}(1+t+t^2)^4\\
y^2 + t^2xy = x^3 + t^2 x             & 12 &  3 & \p{\rm I}_3^* & t^{12}\\
y^2 + t^3xy = x^3 + (1+t)^6x^2        & 12 &  4 & \p{\rm I}_2^* & t^{24}/(1+t)^{12}\\
y^2 + txy = x^3 + (1+t^2)x            &  4 &  1 & {\rm III} & t^8/(1+t^4)\\
y^2 + t^3xy = x^3 + (1+t^6)x          & 12 &  3 & \p{\rm III}^* & t^{24}/(1+t^{12})\\
\hline
\end{array}
$$

For most of these curves, it is easy to determine the behavior of the rational $2$-division point.
As in Section \ref{igusa 3}, one proves:

\begin{proposition}
Let $E_A$ be an elliptic curve over an arbitrary ring $A$ of characteristic two.
Then $E_A$ admits a  section
whose fibers are rational points of order two if and only if it admits
a global Weierstrass equation  of the form $y^2+a_1xy =x^3+a_2x^2+a_4x$
for some   $a_1,a_3,a_4\in A$ with $a_1,a_4$ invertible. The section of order two
is then given by $x=y=0$.
\end{proposition}

\begin{corollary}
\label{2 torsion specialization}
Suppose  $y^2+a_1xy =x^3+a_2x^2+a_4x$ is a minimal Weierstrass equation with nonzero $a_1,a_4\in R$.
Then the rational $2$-division point on $E_K$
has nonzero specialization in the closed fiber of the N\'eron model.
Moreover, we have $a_4\in\maxid_R$ if and only if this
rational point has nonzero class in $\Phi_k$.
\end{corollary}

\begin{theorem}
\mylabel{main result 2}
For all additive Kodaira symbols, there is an elliptic curve $E_K$ containing a rational point of order two
with nonzero specialization in $E_k$ and having the given   reduction type.
For the Kodaira symbols ${\rm III}$, ${\rm III}^*$, ${\rm I}^*_l$, $l\geq 0$, there
are such examples where the specialization has nonzero class in $\Phi_k$, and examples
with zero class in $\Phi_k$. 
\end{theorem}

\proof
For the Kodaira symbols ${\rm II}$, ${\rm II}^*$, ${\rm III}$, ${\rm III}^*$, ${\rm IV}^*$,
the desired examples appear in the table above.
To achieve ${\rm IV}$, consider the elliptic curve
$$
y^2+txy=x^3+a_2x^2+(1+t)x.
$$
This has reduction type ${\rm II}$ for $a_2=0$. For $a_2=t$, the Tate algorithm shows that
the reduction type is ${\rm IV}$, with $\nu(\Delta)=4$ and $\delta=0$.

It remains to treat the cases ${\rm I}^*_l$.
Lemmas \ref{reduction 1} and \ref{reduction 2} easily give the following examples, where $n\geq 2$:
$$
\begin{array}[t]{|l |*{3}{c|}l|}
\hline
\text{minimal Weierstrass equation} &  \text{type} & \text{$2$-torsion}\\
\hline
y^2 + t^{n-1}xy = x^3 + tx^2 + t^nx              &   \p{\rm I}_{2n-4}^*  & (0,0)\\
y^2 + t^{n-1}xy + t^{n-1}(1+t)y= x^3 + t(1+t+t^{n-1})x^2 + t^n(1+t)x       &   \p{\rm I}_{2n-4}^*  & (1+t,(1+t)^2)\\

\hline
y^2 + t^{n-1}xy + t^ny = x^3 + tx^2              &   \p{\rm I}_{2n-3}^*    & (t,0)\\
y^2 + t^nxy + t^n(1+t)y = x^3 + tx^2             &     \p{\rm I}_{2n-3}^*  & (1+t,1+t)\\
\hline
\end{array}
$$
\qed

\section{Semistable reduction in characteristic two}
\mylabel{semistable reduction 2}

Let $R$ be a henselian discrete valuation ring of characteristic $p>0$,
whose residue field $k=R/\maxid_R$ is algebraically closed, and with field of 
fraction $R\subset K$.
Let $E_K$ be an elliptic curve with additive reduction so that $E_K^{(p)}$ has a 
rational $p$-division point.
If $p\geq3$ then we have seen in Theorem \ref{potentially supersingular} and in
Section \ref{igusa 3} that $E_K$ has potentially supersingular reduction.
Although this is not true in characteristic $2$, we see that only additive reduction
of type $I^*_l$ is possible if the curve is not potentially supersingular.

\begin{proposition}
 Let $p=2$ and let $E_K$ be an elliptic curve with additive and potentially ordinary 
 reduction.
 We denote by $\nu(\Delta)$ the valuation of a minimal discriminant and by $\delta$ the
 wild part of the conductor.
 Then there is an integer $d\geq0$ and
 $$
 \begin{array}[t]{|ccc|ccc|}
 \hline
 \multicolumn{3}{|c|}{E_K} & \multicolumn{3}{c|}{E_K^{(2)}}\\
 \nu(\Delta) & \delta & \text{type} & \nu(\Delta) & \delta & \text{type} \\
 \hline
  12d+12 & 4d+2 & I^*_{8d+4} & 
  12d+12 & 4d+2 & I^*_{8d+4} \\
 \hline
 \end{array}
 $$
 and the rational $2$-division point on $E_K^{(2)}$ has nonzero specialization into
 $\Phi_k$.
\end{proposition}

\proof
We may and will assume $R=k[[t]]$.
Since $E_K$ has potentially ordinary reduction,
its $j$-invariant $j=j(E_K)$ lies in $R^\times$ and $E_K$ itself
is ordinary.
In particular, $X_K:y^2+xy=x^3+j^{-1}$
defines an elliptic curve with $j(X_K)=j(E_K)$ and good ordinary
reduction.
Since $E_K$ is ordinary, its automorphism group is $\{\pm1\}$
and hence $X_K$ and $E_K$ differ by a quadratic twist.
As in the proof of Proposition \ref{tautological objects} we conclude 
that $E_K$ is isomorphic to 
$$
y^2+xy=x^3+t^{-2d-1}f^2x^2+j^{-1}
$$
where $f$ is a power series with nonzero constant term.
We have $d\geq0$ because otherwise this curve would have
good reduction.
For $E_K^{(2)}$ we obtain the reduction type, $\nu(\Delta)$ and
$\delta$ analogous to Proposition \ref{pullback types}.
The specialization behavior of the $2$-torsion point
is given by Corollary \ref{2 torsion specialization}.
For $E_K$ the Tate Algorithm shows that 
we have reduction of type $I^*_l$ and $\nu(\Delta)=12d+12$.
Since $\delta(E_K)=\delta(E_K^{(2)})$, we use Ogg's formula
to determine the precise reduction type.
\qed

\medskip
We leave the remaining case to the reader, which is also a generalization
of Proposition \ref{infinity type} and Proposition \ref{infinity pullback}.

\begin{proposition}
 Let $p=2$ and let $E_K$ be an elliptic curve with additive and 
 potentially multiplicative reduction.
 We denote by $\nu(\Delta)$ the valuation of a minimal discriminant and by $\delta$ the
 wild part of the conductor.
 Then there is an integer $d\geq0$ and
 $$
 \begin{array}[t]{|ccc|ccc|}
 \hline
 \multicolumn{3}{|c|}{E_K} & \multicolumn{3}{c|}{E_K^{(2)}}\\
 \nu(\Delta) & \delta & \text{type} & \nu(\Delta) & \delta & \text{type} \\
 \hline
  12d+12-\nu(j) & 4d+2 & I^*_{8d+4-\nu(j)} & 
  12d+12-2\nu(j) & 4d+2 & I^*_{8d+4-2\nu(j)} \\
 \hline
 \end{array}
 $$
 and the rational $2$-division point on $E_K^{(2)}$ has nonzero specialization into
 $\Phi_k$.
\end{proposition}



\begin{thebibliography}{ccccc}





\bibitem{Bosch; Luetkebohmert; Raynaud 1990}
S.\ Bosch, W.\ L\"utkebohmert, M.\ Raynaud:
N\'eron models.
Ergeb.\ Math.\ Grenzgebiete (3) 21.
Springer, Berlin, 1990.


\bibitem{Demazure; Gabriel 1970}
M.\ Demazure, P.\ Gabriel:
Groupes alg\'ebriques. Tome I: G\'eom\'etrie alg\'ebrique, g\'en\'eralit\'es, groupes commutatifs.
North-Holland Publishing Co., 
Amsterdam, 1970. 

\bibitem{Sautoy; Fesenko 2000}
M.\ du Sautoy, I.\ Fesenko:
Where the wild things are: ramification groups and the Nottingham group,
Progr. Math. 184, 287-328, Birkh\"auser, Boston, 2000.

\bibitem{Giraud 1971}
J.\ Giraud:
Cohomologie non ab\'elienne. 
Grundlehren Math.\ Wiss.\ 179. Springer, Berlin, 1971.

\bibitem{Grothendieck 1955}
A.\ Grothendieck:
A general theory of fibre spaces with structure sheaf.
University of Kansas, Department of Mathematics, Report No.\ 4, 1955.

\bibitem{SGA 3b}
M.\ Demazure, A.\ Grothendieck (eds.):
Sch\'emas en groupes II.
Lect.\ Notes Math.\ 152.
Springer, Berlin, 1970.

\bibitem{Fesenko; Vostokov 1993}
I.\ Fesenko,  S.\ Vostokov:
Local fields and their extensions.
Translations of Mathematical Monographs  121. 
American Mathematical Society, Providence, RI, 1993. 

\bibitem{Gunji 1976}
H.\ Gunji:
The Hasse invariant and $p$-division points of an elliptic curve.  
Arch.\ Math.\   27  (1976),   148--158.

\bibitem{Igusa 1968}
J.-I.\ Igusa:
On the algebraic theory of elliptic modular functions.  
J.\ Math.\ Soc.\ Japan  20  (1968) 96--106.

\bibitem{Katz; Mazur 1985}
N.\ Katz, B.\  Mazur:
Arithmetic moduli of elliptic curves.
Annals of Mathematics Studies 108. 
Princeton University Press, Princeton, 1985.



\bibitem{Morandi 1996}
P. Morandi:
Field and Galois theory.
Graduate Texts in Mathematics 167.
Springer, New York, 1996.

\bibitem{Serre 1994}
J.-P.\ Serre:
Cohomologie galoisienne.
Fifth edition. Lect.\ Notes  Math.\ 5. 
Springer, Berlin, 1994.

\bibitem{Shatz 1986}
S.\ Shatz:
Group schemes, formal groups, and p-divisible groups. 
In: 
G.\ Cornell, J.\ Silverman (eds.), Arithmetic geometry, pp.\ 29-78.
Springer, New York, 1986.

\bibitem{Oort; Tate 1970}
J.\ Tate,  F.\ Oort:
Group schemes of prime order.  
Ann.\ Sci.\ \'Ec.\ Norm.\ Sup\'er.\   3, 1--21 (1970).

\bibitem{Raynaud 1974}
M.\ Raynaud:
Sch\'emas en groupes de type $(p,\dots ,p)$. 
Bull.\ Soc.\ Math.\ Fr.\ 102  (1974), 241--280.

\bibitem{Schroeer 2004}
S.\ Schr\"oer:
Some Calabi--Yau threefolds with obstructed deformations
over the Witt vectors.
Compositio Math.\ 140 (2004), 1579--1592.

\bibitem{Serre; Tate 1968}
J.-P.\ Serre, J.\  Tate:
Good reduction of abelian varieties.
Ann.\ Math.\ \ 88 (1968), 492--517.

\bibitem{Serre 1979}
J.-P.\ Serre:
Local fields.
Grad.\ Texts Math.\ 67.
Springer, Berlin, 1979.

\bibitem{Shatz 1970}
S.\ Shatz:.
Finite subschemes of group schemes.
Canad.\ J.\ Math.\ 22 (1970), 1079--1081. 

\bibitem{Silverman 1986}
J.\ Silverman:
The arithmetic of elliptic curves.
Grad.\ Texts Math.\ 106. Springer, Berlin, 1986.

\bibitem{Silverman 1994}
J.\ Silverman:
Advanced topics in the arithmetic of elliptic curves.  
Grad.\ Texts Math.\  151. Springer, New York etc., 1994.

\bibitem{Tate 1975}
J.\ Tate:
Algorithm for determining the type of a singular fiber in an elliptic pencil.  
In: B.\ Birch, W.\ Kuyk (eds.), Modular functions of one variable, IV, pp. 33--52. 
Lecture Notes in Math.\ 476. Springer, Berlin, 1975.

\bibitem{Ulmer 1990a}
D.\ Ulmer:
On universal elliptic curves over Igusa curves.  
Invent.\ Math.\  99  (1990),   377--391.



\end{thebibliography}
\end{document}